\documentclass[11pt]{amsart}

\usepackage{ulem}
\usepackage{graphicx}
\usepackage{amssymb}
\usepackage[top=20mm, bottom=20mm, left=20mm, right=20mm]{geometry}
\usepackage{bm}

\newcommand{\vct}[1]{\bm{\mathsf{#1}}}

\newcommand{\mtx}[1]{\bm{\mathsf{#1}}}

\theoremstyle{definition}
\newtheorem{remark}{Remark}
\numberwithin{definition}{section}

\newcommand{\lsp}{\vspace{3mm}}

\newcommand{\vtwo}[2]{\left[\begin{array}{c} #1 \\ #2 \end{array}\right]}

\newcommand{\mtwo}[4]{\left[\begin{array}{cc}    #1 & #2 \\ #3 & #4  \end{array}\right]}

\begin{document}

\begin{center}
\textbf{Blocked rank-revealing QR factorizations:\\
How randomized sampling can be used to avoid single-vector pivoting}

\lsp

{\small P.G.~Martinsson, Department of Applied Mathematics, University of Colorado at Boulder}

\lsp

May 20, 2015

\lsp

\begin{minipage}{143mm}
\textbf{Abstract:}
Given a matrix $\mtx{A}$ of size $m\times n$, the manuscript describes a algorithm for
computing a QR factorization $\mtx{A}\mtx{P} = \mtx{Q}\mtx{R}$ where $\mtx{P}$ is a
permutation matrix, $\mtx{Q}$ is orthonormal, and $\mtx{R}$ is upper triangular.
The algorithm is blocked, to allow
it to be implemented efficiently. The need for single vector pivoting in classical algorithms
for computing QR factorizations is avoided by the use of randomized sampling to find blocks
of pivot vectors at once. The advantage of blocking becomes particularly pronounced when
$\mtx{A}$ is very large, and possibly stored out-of-core, or on a distributed memory machine.
The manuscript also describes a generalization of the QR factorization where we allow $\mtx{P}$
to be a general orthonormal matrix. In this setting, one can at moderate cost compute a
\textit{rank-revealing} factorization where the mass of $\mtx{R}$ is concentrated to the
diagonal entries. Moreover, the diagonal entries of $\mtx{R}$ closely approximate the
singular values of $\mtx{A}$. The algorithms described have asymptotic flop count
$O(m\,n\,\min(m,n))$, just like classical deterministic methods. The scaling
constant is slightly higher than those of classical techniques, but this is more
than made up for by reduced communication and the ability to block the computation.
\end{minipage}

\end{center}

\lsp

\section{Introduction}
Given an $m\times n$ matrix $\mtx{A}$, with $m \geq n$, the classical QR factorization takes the form
\begin{equation}
\label{eq:AP=QR}
\begin{array}{ccccccccccccc}
\mtx{A} &=& \mtx{Q} & \mtx{R} & \mtx{P}^{*} \\
m\times n && m\times \ell & \ell\times n & n\times n
\end{array}
\end{equation}
where $\mtx{Q}$ is orthonormal, $\mtx{R}$ is upper triangular, and $\mtx{P}$ is a permutation matrix.
The inner dimension $\ell$ can be either $\ell=n$ for an ``economy size'' factorization, or $\ell=m$
for a ``full factorization.''
A standard way of computing the QR factorization is to successively drive $\mtx{A}$ towards upper
triangular form by applying a sequence of Householder reflectors from the left (encoded in $\mtx{Q}$). To ensure that the
diagonal entries of $\mtx{R}$ decay in magnitude, it is common to use column pivoting, which can be
viewed as applying a sequence of permutation matrices representing column swaps to $\mtx{A}$ from
the right (encoded in $\mtx{P}$).

In this manuscript, we use randomized sampling to improve upon the performance of the classical
Householder QR factorization with column pivoting in two ways:
\begin{enumerate}
\item We will enable \textit{blocking} of the algorithm. The resulting algorithm interacts with $\mtx{A}$
via a sequence of BLAS3 operations.
%This makes it easier to optimize the code for high performance when
%$\mtx{A}$ is large, but fits in RAM. It also greatly accelerates computations when $\mtx{A}$ is stored
%``out-of-core'' or on a distributed memory machine.
\item We will describe a variation of the QR factorization for which the diagonal entries of $\mtx{R}$
tend to be \textit{very} close approximations to the singular values of $\mtx{A}$. Such a factorization is commonly
called a \textit{Rank-Revealing QR factorization.}
\end{enumerate}

The computational gains from eliminating column pivoting will be the most pronounced for very large
matrices, in particular ones stored on distributed memory systems, or out-of-core. Likewise, the
RRQR we present is particularly competitive for matrices large enough that computing a full SVD
is not economical. Observe that the methods presented will have \textit{higher} flop counts than
classical methods. The benefit is that they need less data movement.

\section{Notation}
We follow \cite{golub} for general matrix notation. Given a matrix $\mtx{A}$, we let
$\mtx{A}^{*}$ denote its transpose when $\mtx{A}$ is real, and its adjoint when
$\mtx{A}$ is complex. We let $\mtx{I}_{k}$ denote the $k\times k$ identity matrix.

The algorithms presented can advantageously be implemented using standard libraries
for computing matrix factorizations, or matrix-matrix multiplications. This simplifies
coding since it makes the codes portable, and allows us to fully benefit from highly
optimized routines for standard computations. When estimating computational costs, we
use the following simplified model: We assume that the cost of multiplying two matrices
of sizes $\ell\times m$ and $m\times n$ is
$$
C_{\rm mm}\,\ell mn.
$$
The cost of performing a full SVD or QR factorization of a matrix of size $m\times n$, with $m \geq n$, is
$$
C_{\rm qr}\,mn^2,
\qquad\mbox{and}\qquad
C_{\rm svd}\,mn^2.
$$
We will sometimes use non-pivoted QR factorizations. We assume that the cost of this is
$$
C_{\rm qr}^{\rm nopiv}\,mn^2.
$$
Typically, $C_{\rm qr}^{\rm nopiv}$ is smaller that $C_{\rm qr}$.

\begin{remark}
Observe that when computing a pivoted QR factorization of a matrix $\mtx{A}$ of size
$m\times b$ where $m\geq b$, it is always possible to first do a non-pivoted factorization,
and then do a small pivoted factorization on a matrix of size $b\times b$. To be precise,
we first factorize
$$
\begin{array}{cccccccccccccc}
\mtx{A}   &=& \mtx{Q}' & \mtx{R}' \\
m\times b &&  m\times b & b\times b
\end{array}
$$
\textit{with no pivoting.} This is perfectly stable since $\mtx{Q}'$ is built as a product
of ON transforms. Then perform a pivoted QR factorization of the small square matrix $\mtx{R}'$,
$$
\begin{array}{cccccccccccccc}
\mtx{R}' & \mtx{P}   &=& \mtx{Q}'' & \mtx{R}. \\
b\times b & b\times b &&  b\times b & b\times b
\end{array}
$$
Finally, simply set $\mtx{Q} = \mtx{Q}'\mtx{Q}''$ to obtain the factorization (\ref{eq:AP=QR}).
\end{remark}

\section{Review of QR factorization using Householder reflectors}
\label{sec:house}
The algorithm presented in this manuscript is an evolution of the classical technique for computing
a QR factorization via column pivoting and a sequence of Householder reflectors, see, e.g.,
\cite[Sec.~5.2]{golub}. In this section, we briefly review this technique and introduce some notation.
Throughout the section, $\mtx{A}$ is a real matrix of size $m\times n$, with $m\geq n$. The generalization
to complex matrices is trivial.

\subsection{Householder reflectors}
\label{sec:househouse}
Given a vector $\vct{a} \in \mathbb{R}^{k}$, the associated \textit{Householder reflector} $\mtx{H} = \mtx{H}(\vct{a})$ is the
$k\times k$ unitary matrix defined by
$$
\mtx{H} = \mtx{I} - 2\vct{v}\vct{v}^{*},
\qquad\mbox{where}\qquad
\vct{v} = \frac{\beta\,\vct{e}_{1} - \vct{a}}{\|\beta\,\vct{e}_{1} - \vct{a}\|},
\qquad\mbox{and where}\qquad
\beta = -\mbox{sign}(\vct{a}(1))\,\|\vct{a}\|.
$$
The Householder reflector maps $\vct{a}$ to a vector whose entire mass
is concentrated to its first entry:
$$
\mtx{H}\vct{a} = \left[\begin{array}{r}\pm\,\|\vct{a}\| \\ 0 \\ 0 \\ \vdots \\ 0\end{array}\right].
$$

\subsection{QR factorization without pivoting}
\label{sec:housenopivot}
\begin{figure}
\begin{tabular}{cccc}
\includegraphics[height=30mm]{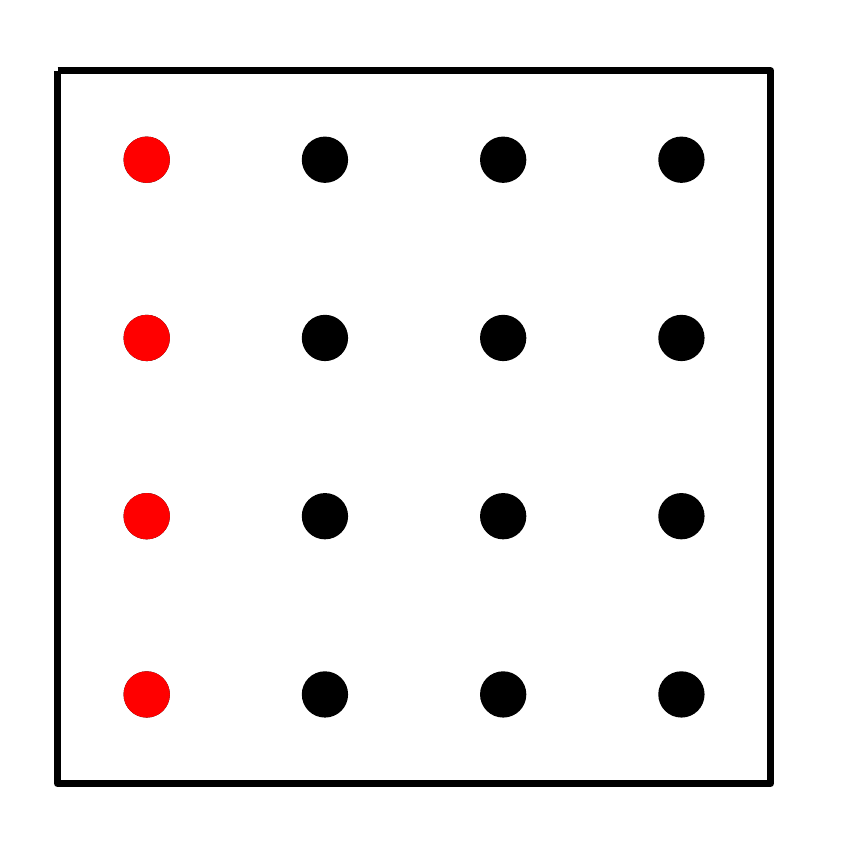} &
\includegraphics[height=30mm]{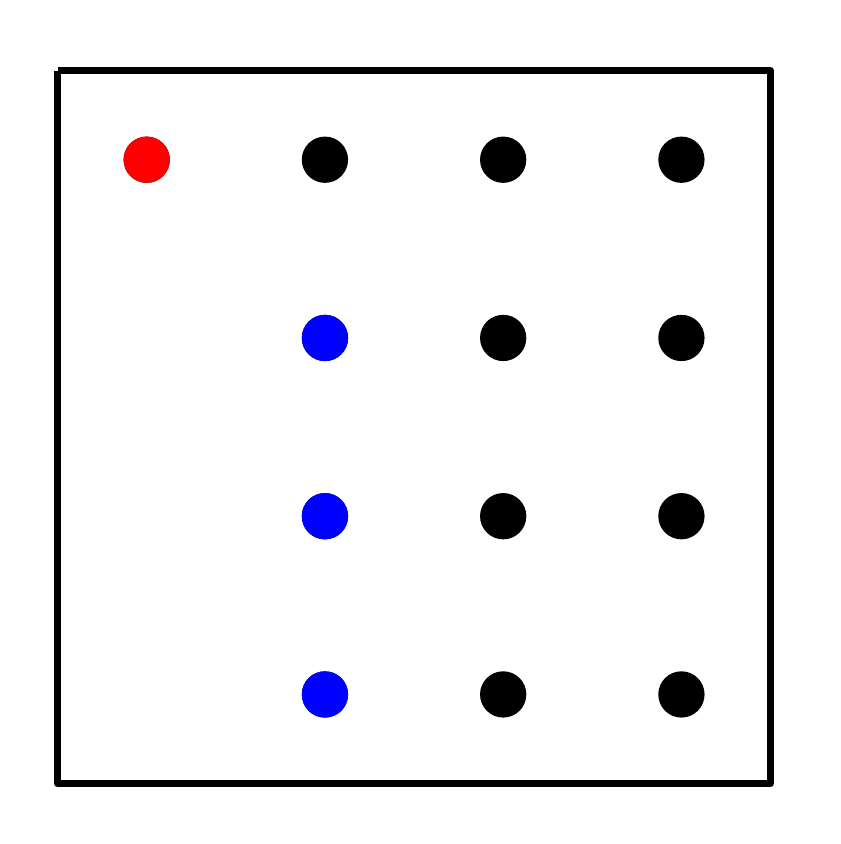} &
\includegraphics[height=30mm]{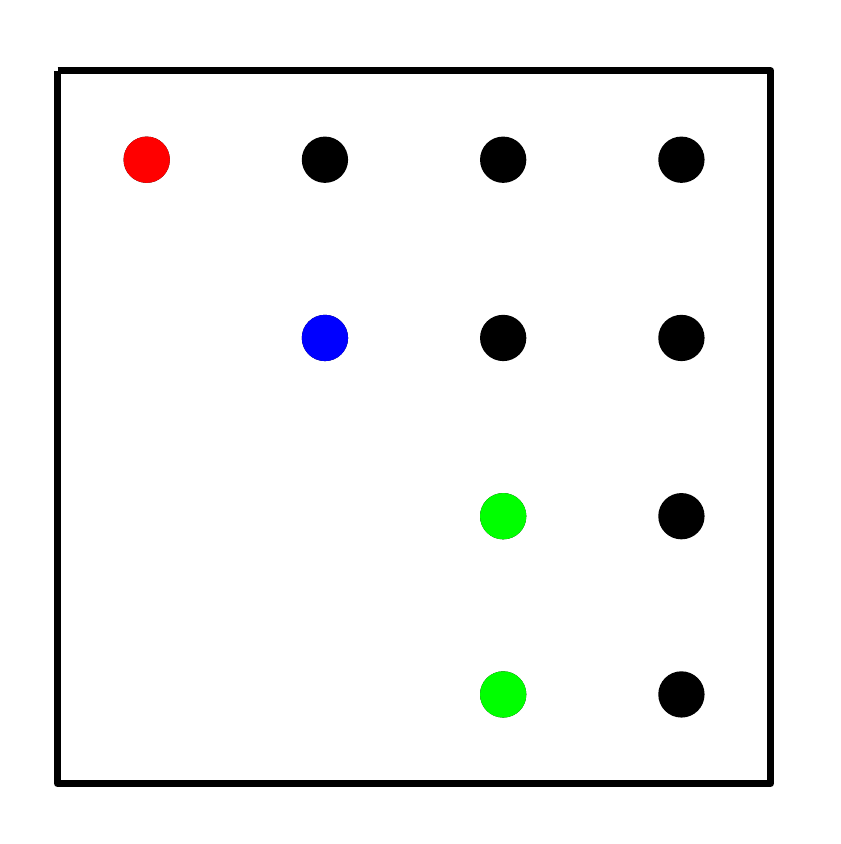} &
\includegraphics[height=30mm]{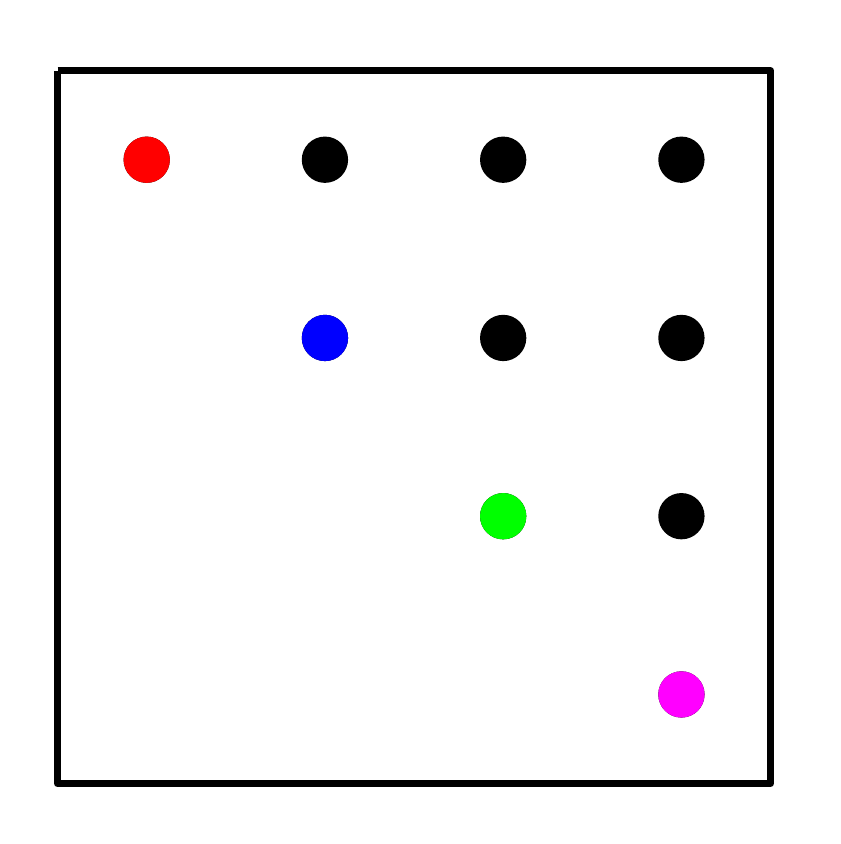} \\
$\mtx{A}_{0}=\mtx{A}$ &
$\mtx{A}_{1}$ &
$\mtx{A}_{2}$ &
$\mtx{A}_{3}=\mtx{R}$
\end{tabular}
\caption{Example of how a $4\times 4$ matrix $\mtx{A}$ is driven to upper triangular form in
a classical QR factorization. Each figure shows the sparsity pattern of the matrix $\mtx{A}_{j}$,
with notation as in Section \ref{sec:house}.}
\label{fig:house}
\end{figure}

In this section, we describe how to drive the given matrix
$\mtx{A}$ to upper triangular form by applying a sequence $n-1$ Householder reflectors from the left.
We set $\mtx{A}_{0} = \mtx{A}$, and let $\mtx{A}_{j}$ denote the result of the first $j$ steps of
the process. The sparsity patterns of these matrices are shown in Figure \ref{fig:house}.
In practice, each $\mtx{A}_{j}$ simply overwrites $\mtx{A}_{j-1}$.

%We will in this section describe a simple process without pivoting that often works, but that can
%fail. The introduction of pivoting in Section \ref{sec:housepivot} will
%make it foolproof and unconditionally stable.

\textit{Step 1:} Let $\mtx{H}$ denote the Householder reflector
associated with the first column of $\mtx{A}_{0}$ (shown in red in Figure \ref{fig:house}).
Set $\mtx{Q}^{(1)} = \mtx{H}$. Then applying $\mtx{Q}^{(1)} = \mtx{Q}^{(1)*}$
to $\mtx{A}_{0}$ will have the effect of ``zeroing out'' out elements below the diagonal in the first column.
We set
\begin{equation}
\label{eq:krux}
\mtx{A}_{1} = \mtx{Q}^{(1)*}\mtx{A}_{0} =
\left[\begin{array}{cccccc}
r_{11} & r_{12}  & r_{13}  & r_{14}  & \cdots \\
0      & a_{22}' & a_{23}' & a_{24}' & \cdots \\
0      & a_{32}' & a_{33}' & a_{34}' & \cdots \\
0      & a_{42}' & a_{43}' & a_{44}' & \cdots \\
\vdots & \vdots  & \vdots  & \vdots
\end{array}\right],
\end{equation}
where $r_{11} = -\mbox{sign}(\mtx{A}_{0}(1,1))\,\|\mtx{A}_{0}(:,1)\|$. (Note: In equation (\ref{eq:krux}), the
matrix $\mtx{Q}^{(1)}$ is symmetric, but we include the transpose to keep notation consistent with
later formulas involving non-symmetric matrices.)

\textit{Step 2:} Let $\mtx{H}$ be the Householder reflector of size $(n-1)\times(n-1)$ associated
with the vector $\mtx{A}_{1}(2:n,2)$ (shown in blue in Figure \ref{fig:house}). Set
$$
\mtx{Q}^{(2)} =
\left[\begin{array}{cc}
\mtx{I}_{1} & \mtx{0} \\
\mtx{0}     & \mtx{H}
\end{array}\right].
$$
The effect of applying
$\mtx{Q}^{(2)*} = \mtx{Q}^{(2)}$ to $\mtx{A}_{1}$ will then be to ``zero out'' all sub-diagonal
elements in the second column. We define
\begin{equation}
\label{eq:krux2}
\mtx{A}_{2} = \mtx{Q}^{(2)*}\mtx{A}_{1} =
\left[\begin{array}{cccccc}
r_{11} & r_{12}  & r_{13}   & r_{14}   & \cdots \\
0      & r_{22}  & r_{23}   & r_{24}   & \cdots \\
0      & 0       & a_{33}'' & a_{34}'' & \cdots \\
0      & 0       & a_{43}'' & a_{44}'' & \cdots \\
\vdots & \vdots  & \vdots   & \vdots
\end{array}\right],
\end{equation}
where $r_{22} = -\mbox{sign}(\mtx{A}_{1}(2,2))\,\|\mtx{A}_{2}(2:n,2)\|$.

\textit{Step $j = 3,\,4,\,\dots,\,n-1$:} Continue just as in Step 2, zeroing out all elements
of $\mtx{A}_{j-1}$ below the diagonal in the $j$'th column, using the Householder reflector
associated with the vector $\mtx{A}_{j-1}(j:m,j)$.

Once all steps have been completed, observe that $\mtx{A}_{n-1}$ is upper triangular, and that we
have now constructed a factorization
$$
\underbrace{\mtx{A}_{n-1}}_{=:\mtx{R}} = \underbrace{\mtx{Q}^{(n-1)*}\,\mtx{Q}^{(n-2)*}\,\cdots\,\mtx{Q}^{(2)*}\,\mtx{Q}^{(1)*}}_{=:\mtx{Q}^{*}}\,\mtx{A}.
$$

\subsection{Column pivoting}
\label{sec:housepivot}
It is often desirable that the diagonal entries of $\mtx{R}$ should form a decreasing sequence
\begin{equation}
\label{eq:Rdecay}
|\mtx{R}(1,1)| \geq |\mtx{R}(2,2)| \geq |\mtx{R}(3,3)| \geq \cdots.
\end{equation}
This can be obtained by introducing \textit{pivoting} into the scheme described in Section \ref{sec:housenopivot}.
The only modification required is that we now need to also hit $\mtx{A}$ with ON-transforms $\{\mtx{P}^{(j)}\}_{j=1}^{n-1}$
from the right.

To be precise, at the start of step $j$, let $j'$ denote the index of the largest column of $\mtx{A}_{j-1}$,
among the ``remaining'' columns $j,\,j+1,\,\dots,n$. Then let $\mtx{P}^{(j)}$ denote the \textit{permutation}
matrix that swaps the columns $j$ and $j'$. Then in the matrix $\mtx{A}_{j-1}\mtx{P}^{(j)}$, the $j$ column will
have the largest column norm among the columns in $\mtx{A}_{j-1}(j:m,j:n)$. Determine the Householder reflector
$\mtx{Q}^{(j)}$ so that it ``zeros out'' the sub-diagonal entries in the $j$'th column of $\mtx{A}_{j-1}\mtx{P}^{(j)}$,
and then set
$$
\mtx{A}_{j} := \mtx{Q}^{(j)*}\,\mtx{A}_{j-1}\,\mtx{P}^{(j)}.
$$

Once all $n-1$ steps have been completed, we end up with a factorization
\begin{equation}
\label{eq:spruce}
\underbrace{\mtx{A}_{n-1}}_{=:\mtx{R}} =
\underbrace{\mtx{Q}^{(n-1)*}\,\mtx{Q}^{(n-2)*}\,\cdots\,\mtx{Q}^{(2)*}\,\mtx{Q}^{(1)*}}_{=:\mtx{Q}^{*}}
\,\mtx{A}\,
\underbrace{\mtx{P}^{(1)}\,\mtx{P}^{(2)}\,\cdots\,\mtx{P}^{(n-2)}\,\mtx{P}^{(n-1)}}_{=:\mtx{P}}.
\end{equation}
Left multiplying (\ref{eq:spruce}) by $\mtx{Q} = \mtx{Q}^{(1)}\mtx{Q}^{(2)}\cdots\mtx{Q}^{(n-1)}$ yields
the factorization (\ref{eq:AP=QR}).

\begin{figure}
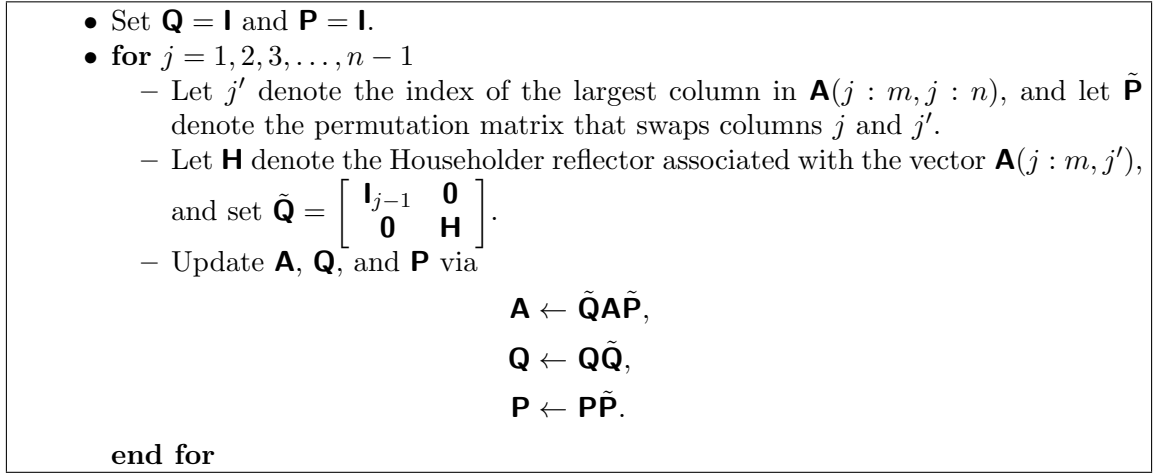

\fbox{\begin{minipage}{150mm}
\begin{itemize}
\item Set $\mtx{Q} = \mtx{I}$ and $\mtx{P} = \mtx{I}$.
\item
\textbf{for} $j = 1,2,3,\dots,n-1$
\begin{itemize}
\item Let $j'$ denote the index of the largest column in $\mtx{A}(j:m,j:n)$, and let
$\tilde{\mtx{P}}$ denote the permutation matrix that swaps columns $j$ and $j'$.
\item Let $\mtx{H}$ denote the Householder reflector associated with the vector $\mtx{A}(j:m,j')$,
and set $\displaystyle \tilde{\mtx{Q}} = \mtwo{\mtx{I}_{j-1}}{\mtx{0}}{\mtx{0}}{\mtx{H}}$.
\item Update $\mtx{A}$, $\mtx{Q}$, and $\mtx{P}$ via
\begin{align*}
\mtx{A}\leftarrow&\ \tilde{\mtx{Q}}\mtx{A}\tilde{\mtx{P}},\\
\mtx{Q}\leftarrow&\ \mtx{Q}\tilde{\mtx{Q}},\\
\mtx{P}\leftarrow&\ \mtx{P}\tilde{\mtx{P}}.
\end{align*}
\end{itemize}
\textbf{end for}
\end{itemize}
\end{minipage}}
\caption{The algorithm \texttt{QR}. Given an input matrix $\mtx{A}$ of size $m\times n$, with $m\geq n$, the
algorithm produces a factorization $\mtx{A}\mtx{P} = \mtx{Q}\mtx{R}$ with $\mtx{R}$ upper triangular,
$\mtx{Q}$ orthonormal, and $\mtx{P}$ a permutation, cf.~(\ref{eq:AP=QR}). The matrix $\mtx{R}$ overwrites $\mtx{A}$.}
\label{fig:basicQR}
\end{figure}

Figure \ref{fig:basicQR} summarizes the classical Householder QR algorithm.

\section{Block pivoting using randomized sampling}
\label{sec:blockpivot}
The QR factorization algorithm based on Householder reflectors described in Section \ref{sec:house}
is exceptionally stable and accurate, but has a shortcoming in that it is hard to \textit{block.}
The traditional column pivoting strategy described necessarily must proceed one vector at a time,
since you cannot find the $j$ pivot column until after the $(j-1)$'th Householder transform has
already been applied. In this section, we address the task of how to find batches of $b$ pivot
vectors at once. Conceptually, we seek to determine a set of $b$ vectors whose spanning volume
is maximal, among the remaining columns.

Suppose we are given a matrix $\mtx{A}$ of size $m\times n$, and let $b$ be a block size.
We will describe two related techniques for computing an ON matrix $\mtx{S}$
of size $n\times n$ such that the first $b$ columns of $\mtx{A}\mtx{S}$ form good choices
for the first $b$ ``pivot columns'' in a QR factorization. In Section \ref{sec:pivotperm},
we describe a technique for constructing a matrix $\mtx{S}$ that is a permutation matrix, so
that the first $k$ columns of $\mtx{A}\mtx{S}$ are simply $k$ chosen columns from $\mtx{A}$.
In Section \ref{sec:pivotON} we generalize slightly further from the classical QR factorization,
and will describe a ``pivoting matrix'' $\mtx{S}$ that is ON, but is not merely a permutation
(in fact, it will consist of a sequence of Householder reflectors). In other words, each of
the first $k$ columns of $\mtx{A}\mtx{S}$ will consist of \textit{linear combinations} of
columns of $\mtx{A}$.

\subsection{Finding a permutation matrix $\mtx{S}$}
\label{sec:pivotperm}
Draw a Gaussian random matrix $\mtx{\Omega}$ of size $m\times b$, and compute a
``sampling matrix'' $\mtx{Y}$ via
\begin{equation}
\label{eq:defY}
\begin{array}{cccccccccc}
\mtx{Y} &=& \mtx{A}^{*}&\mtx{\Omega}.\\
n\times b && n\times m & m\times b
\end{array}
\end{equation}
Then one can prove that with very high probability, the linear dependencies among
the columns of $\mtx{Y}$ closely match the linear dependencies among the columns of
$\mtx{A}$. This means that if we perform a classical QR factorization
\begin{equation}
\label{eq:Yrows}
\begin{array}{cccccccccccccccc}
\mtx{Y}^{*}&\mtx{S} &=& \mtx{Q}&\mtx{R},\\
b \times n & n\times n && b \times b & b\times n
\end{array}
\end{equation}
then the permutation matrix $\mtx{S}$ chosen will also be a good permutation
matrix for $\mtx{A}$. The cost of finding this matrix $\mtx{S}$ is
$$
C_{\rm mm}mnb + C_{\rm qr}nb^{2}.
$$

\subsection{Finding a general orthonormal matrix $\mtx{S}$}
\label{sec:pivotON}
Construct a Gaussian random matrix $\mtx{\Omega}$ of size $m\times b$, and again
compute a ``sampling matrix'' $\mtx{Y} = \mtx{A}^{*}\mtx{\Omega}$, cf.~(\ref{eq:defY}).
Next, perform a QR factorization of $\mtx{Y}$ to form the factorization
\begin{equation}
\label{eq:Ycols}
\begin{array}{cccccccccccccccc}
\mtx{Y}&\mtx{P} &=& \mtx{S}&\mtx{R}.\\
n \times b & b\times b && n \times n & n\times b
\end{array}
\end{equation}
Comparing (\ref{eq:Ycols}) to (\ref{eq:Yrows}) we observe that (\ref{eq:Ycols}) represents
a orthonormalization of the \textit{columns} of $\mtx{Y}$, while (\ref{eq:Yrows}) orthonormalizes
the \textit{rows} of $\mtx{Y}$. Now observe that
$$
\mtx{S}^{*}\mtx{Y} = \left[\begin{array}{c}\mtx{R} \\ \mtx{0}\end{array}\right].
$$
In other words, $\mtx{S}^{*}$ is on orthonormal map that rotates all the mass in $\mtx{Y}$ into
the first $k$ rows. This means that in the matrix $\mtx{A}\mtx{S}$, the leading $b$ columns
approximately span the same space as the leading $b$ left singular vectors of $\mtx{A}$, which
would form the ``ideal'' pivoting vectors.

\begin{remark}
One can prove that in computing the factorization (\ref{eq:Ycols}), it is possible to forgo
pivoting entirely (so that $\mtx{P} = \mtx{I}$), which will accelerate this computation.
\end{remark}

\begin{remark}
The matrix $\mtx{S}$ described in this subsection is a large (of size $n\times n$) ON matrix.
Note, however, that it is composed simply of a product of $b$ Householder reflectors. This means
that the factorization (\ref{eq:Ycols}) can be computed in $O(nb^{2})$ operations, $\mtx{S}$ requires
$O(nb)$ storage, and can be applied to a vector using $O(nb)$ flops.
\end{remark}

\subsection{Over-sampling}
\label{sec:over}
The accuracy of the procedures described in this section can be improved
by constructing a few ``extra samples.'' Let $r$ denote an over-sampling parameter. The choice
$r=10$ is often excellent, and $r=b$ leads to very high accuracy. Then generate a Gaussian
matrix $\mtx{\Omega}$ of size $m\times (b+r)$. Then, in computing the factorizations
(\ref{eq:Yrows}) or (\ref{eq:Ycols}), execute only the first $b$ steps of the QR process.

\section{Blocked QR}
\label{sec:blockQR}

We describe a blocked version of the basic Householder QR algorithm (cf.~Section \ref{sec:house})
in Section \ref{sec:blockQRsub}. The block pivoting can be done using either permutation matrices
(cf.~Section \ref{sec:pivotperm}) or Householder reflectors (cf.~Section \ref{sec:pivotON}). The
two resulting algorithms are analyzing in Sections \ref{sec:blockperm} and \ref{sec:blockON}, respectively.

\subsection{The algorithm}
\label{sec:blockQRsub}
\begin{figure}
\begin{tabular}{ccccc}
\includegraphics[height=32mm]{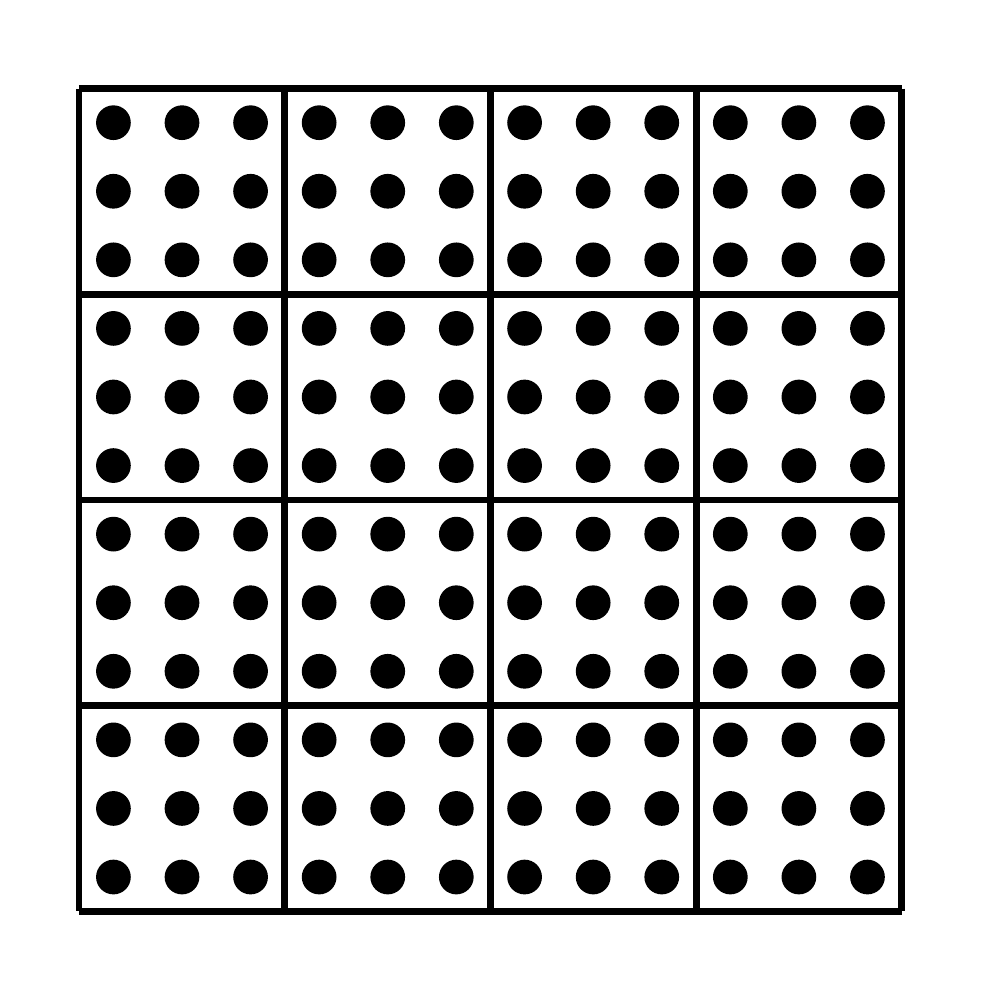} &
\includegraphics[height=32mm]{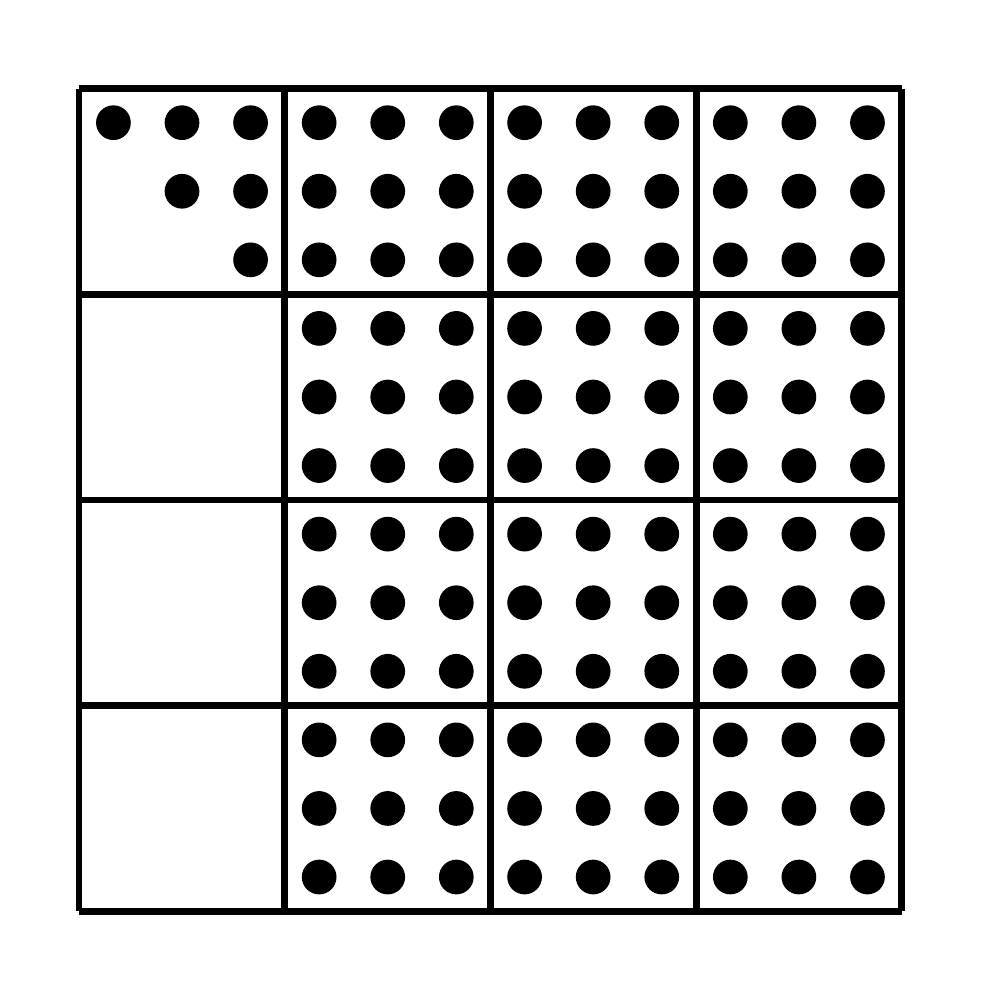} &
\includegraphics[height=32mm]{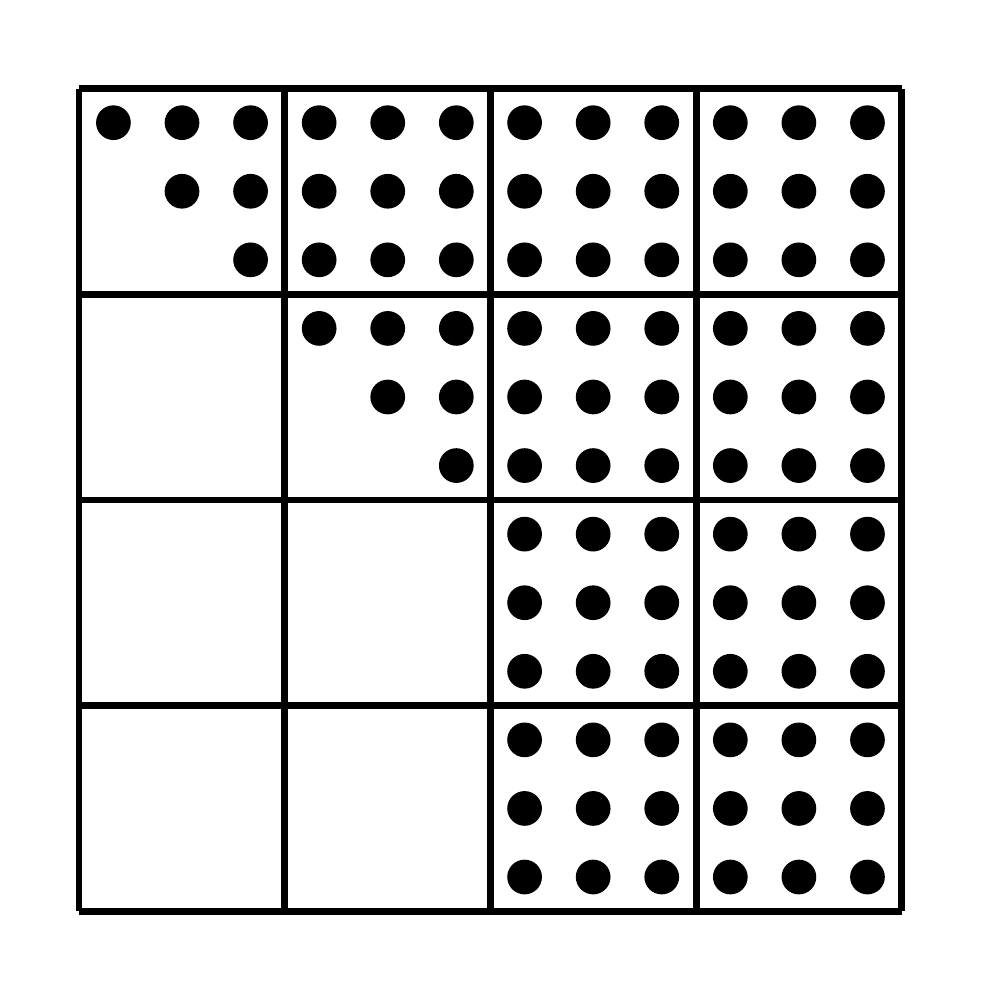} &
\includegraphics[height=32mm]{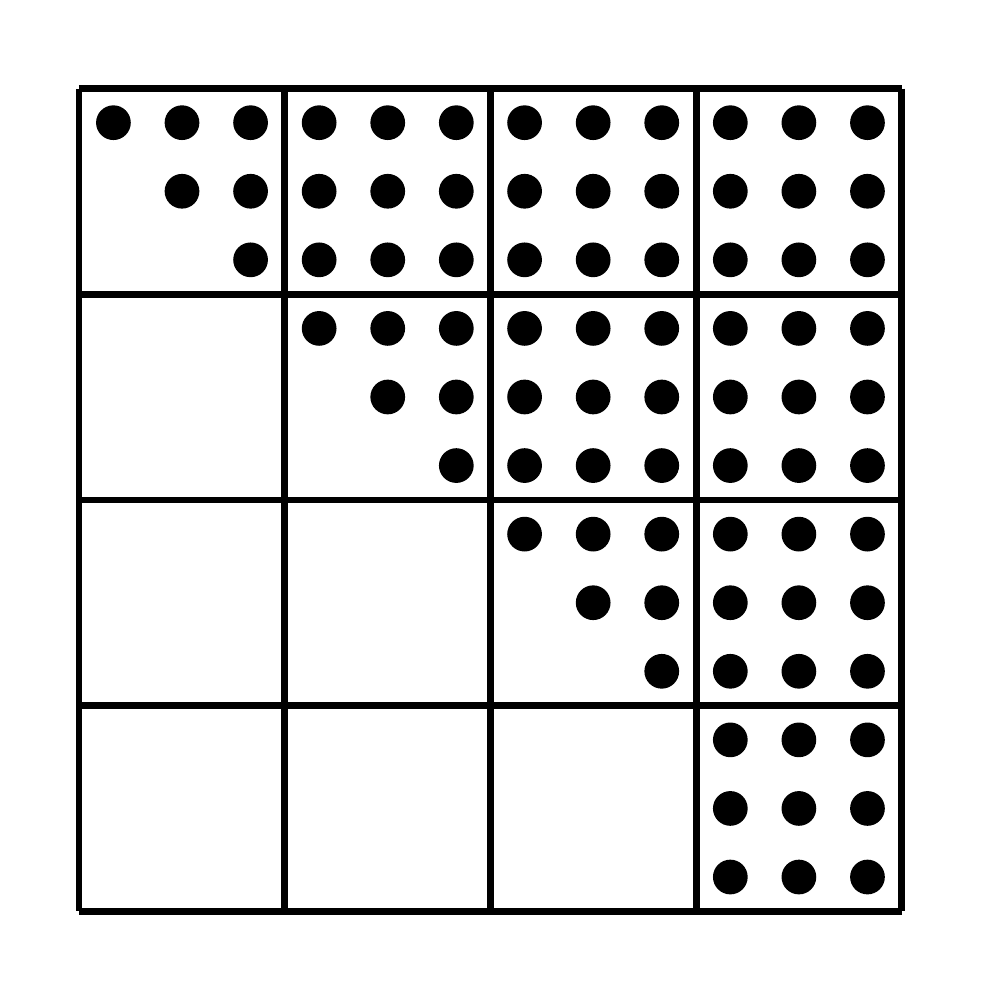} &
\includegraphics[height=32mm]{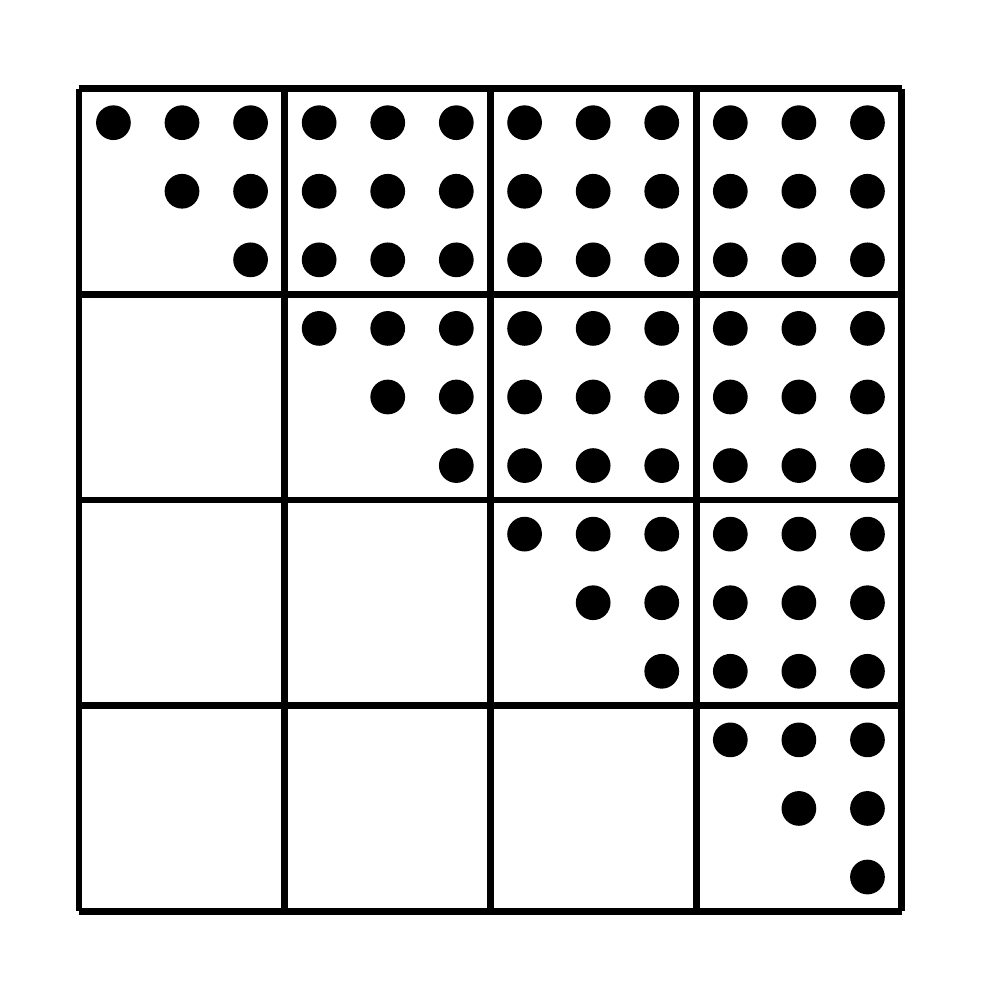} \\
$\mtx{A}_{0}=\mtx{A}$ &
$\mtx{A}_{1}$ &
$\mtx{A}_{2}$ &
$\mtx{A}_{3}$ &
$\mtx{A}_{4}=\mtx{R}$
\end{tabular}
\caption{Example of how a $12\times 12$ matrix $\mtx{A}$ is driven to upper triangular form in
the blocked randomized QR factorization. Each figure shows the sparsity pattern of matrix $\mtx{A}_{j}$,
with notation as in Section \ref{sec:blockQR}.}
\label{fig:blockQRsparsity}
\end{figure}

Suppose that $\mtx{A}$ is an $m\times n$ matrix, with $m\geq n$. We will drive $\mtx{A}$ to
upper triangular form by processing blocks of $b$ vectors at a time. Suppose for simplicity that
$n$ is a multiple of $b$, so that $n = bp$ for some integer $p$. The blocked algorithm
resulting is shown in detail in Figure \ref{fig:blockQR}. The non-zero elements in the matrix
$\mtx{A}_{j}$ obtained after $j$ steps of the blocked algorithm is shown in Figure \ref{fig:blockQRsparsity}.

\begin{figure}
\fbox{\begin{minipage}{150mm}
\begin{itemize}
\item Set $\mtx{Q} = \mtx{I}$ and $\mtx{P} = \mtx{I}$.
\item
\textbf{for} $i = 1,2,3,\dots,p-1$
\begin{itemize}
\item Partition the index vector $I = [I_{1}\ I_{2}\ I_{3}]$ so that $I_{2} = b(i-1) + (1:b)$ is the block processed in step $i$.
Then partition the matrix accordingly,
$$
\mtx{A} =
\left[\begin{array}{ccc}
\mtx{A}_{11} & \mtx{A}_{12} & \mtx{A}_{13} \\
\mtx{0}      & \mtx{A}_{22} & \mtx{A}_{23} \\
\mtx{0}      & \mtx{A}_{32} & \mtx{A}_{33}
\end{array}\right].
$$
(Observe that in the first step, $I_{1} = []$ and $\mtx{A}$ has only $2\times 2$ blocks.)
\item Determine a pivoting matrix $\tilde{\mtx{P}}$ by processing $\mtx{A}([I_{2}\ I_{3}],[I_{2}\ I_{3}])$,
as described in Section \ref{sec:blockpivot}.
Then update the last two block columns of $\mtx{A}$ accordingly,
$$
\left[\begin{array}{ccc}
\mtx{A}_{12}' & \mtx{A}_{13}' \\
\mtx{A}_{22}' & \mtx{A}_{23}' \\
\mtx{A}_{32}' & \mtx{A}_{33}'
\end{array}\right] =
\left[\begin{array}{ccc}
\mtx{A}_{12} & \mtx{A}_{13} \\
\mtx{A}_{22} & \mtx{A}_{23} \\
\mtx{A}_{32} & \mtx{A}_{33}
\end{array}\right]\,\tilde{\mtx{P}}.
$$
\item Execute a QR-factorization
$$
\begin{array}{ccccccccccc}
\vtwo{\mtx{A}_{22}'}{\mtx{A}_{23}'} & \hat{\mtx{P}} &=& \tilde{\mtx{Q}} & \vtwo{\mtx{R}_{22}}{\mtx{0}},\\
m'\times b & b\times b && m' \times m' & m'\times b
\end{array}
$$
where $m' = m - (i-1)b$ is the length of the index vector $[I_{2}\ I_{3}]$. Observe that while $\tilde{\mtx{Q}}$ is
large, it consists simply of a product of $b$ Householder reflectors.
\item Compute the new blocks $\mtx{A}_{23}''$ and $\mtx{A}_{33}''$ via
$\displaystyle
\vtwo{\mtx{A}_{23}''}{\mtx{A}_{33}''} = \tilde{\mtx{Q}}^{*}
\vtwo{\mtx{A}_{23}'}{\mtx{A}_{33}'}.
$
\item Update the matrices $\mtx{A}$, $\mtx{Q}$, and $\mtx{P}$, via
\begin{align*}
\mtx{A}\leftarrow&\
\left[\begin{array}{ccc}
\mtx{A}_{11} & \mtx{A}_{12}'\hat{\mtx{P}} & \mtx{A}_{13}'  \\
\mtx{0}      & \mtx{R}_{22}             & \mtx{A}_{23}'' \\
\mtx{0}      & \mtx{0}                  & \mtx{A}_{33}''
\end{array}\right],\\
\mtx{Q}\leftarrow&\ \mtx{Q}\left[\begin{array}{cc} \mtx{I}_{n_1} & \mtx{0} \\ \mtx{0} & \tilde{\mtx{Q}}\end{array}\right],\\
\mtx{P}\leftarrow&\ \mtx{P}\left[\begin{array}{cc} \mtx{I}_{n_1} & \mtx{0} \\ \mtx{0} & \tilde{\mtx{P}}\end{array}\right]\,
\left[\begin{array}{ccc}
\mtx{I}_{n_1} & \mtx{0}       & \mtx{0} \\
\mtx{0}     & \hat{\mtx{P}} & \mtx{0} \\
\mtx{0}     & \mtx{0}       & \mtx{I}_{n_3}
\end{array}\right].
\end{align*}
\end{itemize}
\textbf{end for}
\item At this point, all that remains is to process the lower right $b\times b$ block. Partition
$$
\mtx{A} = \mtwo{\mtx{A}_{11}}{\mtx{A}_{12}}{\mtx{0}}{\mtx{A}_{22}}
$$
so that $\mtx{A}_{22}$
is of size $b\times b$. Observe that $\mtx{A}_{11}$ is already upper triangular. Now compute a QR factorization
$\tilde{\mtx{Q}}\mtx{R}_{22} = \mtx{A}_{22}\tilde{\mtx{P}}$. Then simply update
$$
\mtx{A}\leftarrow\mtwo{\mtx{A}_{11}}{\mtx{A}_{12}}{\mtx{0}}{\mtx{R}_{22}},\qquad
\mtx{Q} = \mtx{Q}\mtwo{\mtx{I}_{n_1}}{\mtx{0}}{\mtx{0}}{\tilde{\mtx{Q}}},\qquad
\mtx{P} = \mtx{P}\mtwo{\mtx{I}_{n_1}}{\mtx{0}}{\mtx{0}}{\tilde{\mtx{P}}}.
$$
\end{itemize}
\end{minipage}}
\caption{The algorithm \texttt{blockQR}. Given an input matrix $\mtx{A}$ of size $m\times n$, with $m\geq n$, the
algorithm produces a factorization $\mtx{A}\mtx{P} = \mtx{Q}\mtx{R}$ with $\mtx{R}$ upper triangular,
and $\mtx{P}$ and $\mtx{Q}$ orthonormal, cf.~(\ref{eq:AP=QR}). The matrix $\mtx{R}$ overwrites $\mtx{A}$.}
\label{fig:blockQR}
\end{figure}

For efficiency, it is \textit{essential} that all manipulations involving ON matrices be executed
using the fact that these are products of Householder reflectors. Observe that if an $n\times n$ matrix
$\mtx{U}$ is a product of $b$ Householder reflectors, then $\mtx{U}$ admits the representation
$$
\begin{array}{ccccccccccccc}
\mtx{U} &=& \mtx{I} &+&\mtx{W}&\mtx{Y}\\
n\times n && n\times n && n\times b & b\times n
\end{array}
$$
for some matrices $\mtx{W}$ and $\mtx{Y}$, see \cite{1987_bischof_WY}. This should be exploited whenever
a matrix such as $\mtx{U}$ is applied to $\mtx{A}$.

\subsection{Block pivoting via permutation matrices}
\label{sec:blockperm}
Suppose that we build the pivoting matrix $\tilde{\mtx{P}}$ in \texttt{blockQR} as a
permutation matrix, using the strategy of Section \ref{sec:pivotperm}. Observe that
this step still requires ``single-vector'' column pivoting. The gain here is that the
pivoted factorization is executed on a small matrices of size $b\times n$. A secondary
pivoted factorization is also executed to build the matrices $\tilde{\mtx{Q}}$ and $\hat{\mtx{P}}$,
but again, the matrix involved is of size $(m-(i-1)b)\times b$.

The key point is that the interaction with $\mtx{A}$ is done exclusively via matrix-matrix
multiplications. This leads to a modest acceleration when $\mtx{A}$ fits in RAM but is large.
When $\mtx{A}$ is stored on a distributed memory machine, or is stored out-of-core, we expect
a decisive speed-up.

\subsection{Block pivoting via Householder reflectors}
\label{sec:blockON}
The asymptotic cost of the algorithm using Householder reflectors instead of permutation matrices
has the same scaling as the algorithm based on permutation matrices, although the constants are
slightly larger. The benefit of this algorithm is that the pivoting is ``better'' in the sense
that more of the mass is concentrated to the diagonal blocks, as we will discuss in detail in
Section \ref{sec:RRQR}.

\section{Rank-Revealing QR factorization}
\label{sec:RRQR}

The idea of \textit{Rank-Revealing QR factorization (RRQR)} \cite{gu1996,1987_chan_RRQR}
is to combine some of the
advantages of QR factorizations and singular value decompositions. While the SVD is
excellent for revealing how well a given matrix can be approximated by a matrix of
low rank, it can in principle only be computed via iterative procedures. In practice,
the best iterative schemes converge fast enough that they in \textit{most} circumstances
behave just like deterministic algorithms. The idea of an RRQR is that it can be computed
using non-iterative methods (typically substantially faster than a full SVD) and reveals
rank almost as well as an SVD.

%The requirement that the diagonal entries of $\mtx{R}$ should decay, cf.~(\ref{eq:Rdecay}), is often
%sufficient. However, it can be substantially improved, to the point where the diagonal entries to
%some degree approximate the singular values of $\mtx{A}$, see, e.g., \cite{gu1996,1987_chan_RRQR}.

\subsection{Low rank approximation}
Let $\mtx{A}$ be of size $m\times n$, with $m \geq n$. The Singular Value Decomposition (SVD)
of $\mtx{A}$ takes the form
\begin{equation}
\label{eq:SVD}
\begin{array}{ccccccccccc}
\mtx{A} &=& \mtx{U} & \mtx{D} & \mtx{V}^{*},\\
m\times n && m\times n & n\times n & n\times n
\end{array}
\end{equation}
where $\mtx{U}$ and $\mtx{V}$ are orthonormal. The diagonal matrix $\mtx{D}$ has as its
diagonal entries the singular values $\{\sigma_{j}\}_{j=1}^{n}$, ordered so that
$\sigma_{1} \geq \sigma_{2} \geq \cdots \geq \sigma_{n} \geq 0$. The Eckart-Young
theorem \cite{1936_eckart_young} states that when matrices are measured in either the spectral
or the Frobenius norm, then the truncated SVD is an optimal low-rank approximation to a matrix.
To be precise, fix a rank $k$, and partition
\begin{equation}
\label{eq:SVDpart}
\mtx{A} = \bigl[\mtx{U}_{1}\ \mtx{U}_{2}\bigr]\,
\mtwo{\mtx{D}_{11}}{\mtx{0}}{\mtx{0}}{\mtx{D}_{22}}\,
\vtwo{\mtx{V}_{1}^{*}}{\mtx{V}_{2}^{*}}
\end{equation}
so that $\mtx{D}_{11}$ is of size $k\times k$. Then in the spectral norm, we have
$$
\min\{\|\mtx{A} - \mtx{B}\|\,:\,\mtx{B}\mbox{ has rank }k\} =
\|\mtx{A} - \mtx{U}_{1}\mtx{D}_{11}\mtx{V}_{1}^{*}\| =
\|\mtx{U}_{2}\mtx{D}_{22}\mtx{V}_{2}^{*}\| =
\|\mtx{D}_{22}\| = \sigma_{k+1}.
$$
The statement for the Frobenius norm is analogous:
$$
\min\{\|\mtx{A} - \mtx{B}\|_{\rm Fro}\,:\,\mtx{B}\mbox{ has rank }k\} =
\|\mtx{A} - \mtx{U}_{1}\mtx{D}_{11}\mtx{V}_{1}^{*}\|_{\rm Fro} =
\|\mtx{U}_{2}\mtx{D}_{22}\mtx{V}_{2}^{*}\|_{\rm Fro} =
\|\mtx{D}_{22}\|_{\rm Fro} = \left(\sum_{j=k+1}^{n}\sigma_{j}^{2}\right)^{1/2}.
$$

Now suppose that we do the same thing for a QR factorization
$$
\mtx{A}\mtx{P} =
\bigl[\mtx{Q}_{1}\ \mtx{Q}_{2}\bigr]\,
\mtwo{\mtx{R}_{11}}{\mtx{R}_{12}}{\mtx{0}}{\mtx{R}_{22}}.
$$
Then the error in a rank-$k$ approximation would be
$$
\|\mtx{A}\mtx{P} - \mtx{Q}_{1}\bigl[\mtx{R}_{11}\ \mtx{R}_{12}\bigr]\| =
\|\mtx{Q}_{2}\mtx{R}_{22}\| =
\|\mtx{R}_{22}\|.
$$
In order for the error in the truncated QR to be close to optimal, we
would like to have that $\|\mtx{R}_{22}\| \approx \|\mtx{D}_{22}\|$.
Enforcing this for every choice of $k$, we find that we seek
\begin{equation}
\label{eq:R22_cond}
\|\mtx{R}((k+1):n,(k+1):n)\| \approx \|\mtx{D}((k+1):n,(k+1):n)\|,
\qquad k = 1,\,2,\,\dots, n-1.
\end{equation}
A closely related condition is the slightly stronger, and simpler,
condition that the diagonal entries of $\mtx{R}$ should all
approximate the corresponding singular values, so that
\begin{equation}
\label{eq:Dkk_cond}
|\mtx{R}(k,k)| \approx \sigma_{k},
\qquad k = 1,\,2,\,\dots, n.
\end{equation}
Informally, the idea is to move as much mass as possible in the
matrix $\mtx{R}$ onto the diagonal entries.

In traditional QR factorizations, one typically assumes that the
matrix $\mtx{P}$ is a permutation matrix, so that the columns of
$\mtx{Q}(:,1:k)$ provide an orthonormal basis for a selection of
$k$ columns of $\mtx{A}$. Under this condition, it is possible to
construct counter-examples that demonstrate that the best possible
QR factorization will fail to achieve either (\ref{eq:R22_cond})
or (\ref{eq:Dkk_cond}). However, since we allow $\mtx{P}$ to be a
general orthonormal matrix, there is nothing in principle that
prevents us from realizing these bounds to very high precision.
We will in Sections \ref{sec:diagonalizing} and \ref{sec:power}
describe three modifications to the QR factorization scheme in
Section \ref{sec:blockQR} that will make the scheme output a
very high quality QR factorization satisfying (\ref{eq:R22_cond})
and (\ref{eq:Dkk_cond}). The final scheme uses the following
building blocks (where $b$ is the block size):
\begin{itemize}
\item The matrix $\mtx{A}$ will be interacted with only via
matrix-matrix multiplies involving matrices with at most $2b$
columns or rows, and low-rank updates.
\item We will use \textit{unpivoted} QR factorizations of
matrices of size at most $m\times 2b$ or $n\times 2b$.
\item We will compute full SVDs of $n/b$ matrices of size at most $2b\times 2b$.
\end{itemize}

\subsection{Diagonalizing the diagonal blocks}
\label{sec:diagonalizing}

In the blocked QR algorithm in Figure \ref{fig:blockQR}, we can at very low
cost enforce that the diagonal block $\mtx{R}_{22}$ be not only upper triangular,
but \textit{diagonal.} All that is required is to replace the local QR factorization
by a (full) SVD
\begin{equation}
\label{eq:localSVD}
\begin{array}{cccccccccccc}
\vtwo{\mtx{A}_{22}'}{\mtx{A}_{32}'} &=& \tilde{\mtx{U}} & \vtwo{\mtx{D}_{22}}{\mtx{0}} & \tilde{\mtx{V}}^{*}, \\
m' \times b && m' \times m' & m' \times b & b\times b
\end{array}
\end{equation}
and then use $\tilde{\mtx{U}}$ instead of $\tilde{\mtx{Q}}$, and $\tilde{\mtx{V}}$ instead of $\hat{\mtx{P}}$.
We compute the factorization (\ref{eq:localSVD}) via two steps: First, perform an unpivoted QR factorization
$$
\begin{array}{cccccccccccc}
\vtwo{\mtx{A}_{22}'}{\mtx{A}_{32}'} &=& \tilde{\mtx{Q}} & \vtwo{\mtx{R}_{22}}{\mtx{0}}. \\
m' \times b && m' \times m' & m' \times b
\end{array}
$$
Observe that $\tilde{\mtx{Q}}$ is a product of $b-1$ Householder reflectors. Then compute the SVD of $\mtx{R}_{22}$:
$$
\begin{array}{cccccccccccc}
\mtx{R}_{22} &=& \mtx{U}' & \mtx{D}_{22} & \tilde{\mtx{V}}^{*}. \\
b \times b && b \times b & b \times b & b\times b
\end{array}
$$
Finally, compute $\tilde{\mtx{U}}$ via
$$
\tilde{\mtx{U}} =
\tilde{\mtx{Q}}\,
\mtwo{\mtx{U}'}{\mtx{0}}{\mtx{0}}{\mtx{I}_{n_3}}.
$$
The total cost of this step is
$$
C_{\rm qr}^{\rm nopiv}\,m'\,b^{2} + C_{\rm svd}\,b^{3} + C_{\rm mm}\,m'\,b^{2}.
$$

\subsection{Power iteration to improve pivoting further}
\label{sec:power}
In Section \ref{sec:pivotON}, we describe a technique for finding an ON matrix $\mtx{S}$
such that the first $b$ columns of $\mtx{A}\mtx{S}$ form suitable ``pivots'' in a
QR factorization. A theoretically excellent choice for such a matrix $\mtx{S}$ would be
an ON matrix whose first $b$ columns span the leading $b$ right singular vectors of $\mtx{A}$.
To see why, suppose that we partition $\mtx{S}$ so that
$$
\mtx{S} = \bigl[\mtx{S}_{1}\ \mtx{S}_{2}\bigr]
$$
in such a way that $\mtx{S}_{1}$ is of size $n\times b$, and
\begin{equation}
\label{eq:VS}
\mtx{V}_{1}^{*}\mtx{S}_{2} = \mtx{0},
\qquad\mbox{and}\qquad
\mtx{V}_{2}^{*}\mtx{S}_{1} = \mtx{0}.
\end{equation}
With the SVD of $\mtx{A}$ partitioned as in (\ref{eq:SVDpart}), we then find that
$$
\mtx{A}\mtx{S} =
\bigl[\mtx{U}_{1}\mtx{D}_{11}\mtx{V}_{1}^{*}\mtx{S}_{1}\ \
      \mtx{U}_{2}\mtx{D}_{22}\mtx{V}_{2}^{*}\mtx{S}_{1}\bigr].
$$
Then the pivot columns will span precisely the leading $b$ left singular vectors,
which means that after we apply the first $b$ Householder reflectors from the left,
the resulting matrix will be block diagonal.

Now, finding a matrix $\mtx{S}$ for which (\ref{eq:VS}) holds is hard, since
it amounts to finding the span of the leading $b$ right singular vectors of
$\mtx{A}$. (The purpose of computing an RRQR is precisely to avoid this!)
But finding an \textit{approximate} span of the $b$ leading singular vectors
is something that randomized sampling excels at. The matrix $\mtx{Y}$
described in Section \ref{sec:pivotON} is built specifically so that its
columns span the space we seek to determine. The alignment between the range
of $\mtx{Y}$ and the range of $\mtx{V}_{1}$ can be further improved by applying
a power of the remaining columns. To be precise, let $\mtx{X} = \mtx{A}([I_2,I_3],[I_2,I_3])$
denote the block of $\mtx{A}$ that remains to be driven to upper triangular form.
Then if we fix an integer $q$, and build a sample matrix
$$
\mtx{Y} = \bigl(\mtx{X}^{*}\mtx{X}\bigr)\,\mtx{X}^{*}\,\mtx{\Omega},
$$
where $\mtx{\Omega}$ is again a matrix with i.i.d.~Gaussian entries, then as $q$
increases, the range of $\mtx{Y}$ tends to rapidly converge to the range of
$\mtx{V}_{1}$ (e.g.~if all singular values are distinct, then such convergence
can easily be proven). In practice, choosing $q=1$ or $q=2$ tends to give
excellent results.

Finally, to attain truly high accuracy, we will employ over-sampling as described
in Section \ref{sec:over}, but more aggressively than before. While in the standard
QR factorization, it is fine to choose the over-sampling parameter $p$ to be a fixed
small integer (say $p=5$ or $p=10$, or even $p=0$), we have empirically found that
with $p=b$, we attain an excellent alignment between the spans of $\mtx{V}_{1}$ and
$\mtx{S}_{1}$.

%\subsection{Background on RRQRs}
%The singular value decomposition (SVD) is in many ways an ``ideal'' matrix factorization
%for determining how accurately a given matrix can be approximated by a matrix of low rank.
%The idea of RRQRs is that while the SVD necessarily requires an iterative algorithm for
%its computation, an RRQR can be computed using deterministic techniques. Moreover, the
%while the cost of computing an RRQR and an SVD scale in the same manner, in practice, the
%RRQR is typically cheaper to compute. To illustrate the questions that arise, let us
%compare the SVD and the QR factorization for an $m\times n$ matrix $\mtx{A}$ (with
%$m \geq n$, as usual):
%\begin{equation}
%\label{eq:compare}
%\begin{array}{ccccccccccc}
%\mtx{A} &=& \mtx{Q} & \mtx{R} & \mtx{P}^{*}\\
%m\times n && m\times n & n\times n & n\times n
%\end{array}
%\qquad\mbox{versus}\qquad
%\begin{array}{ccccccccccc}
%\mtx{A} &=& \mtx{U} & \mtx{D} & \mtx{V}^{*}\\
%m\times n && m\times n & n\times n & n\times n
%\end{array}
%\end{equation}
%In the two factorizations ...

\subsection{The algorithm \texttt{blockRRQR}}

Our method for computing an RRQR is a obtained by starting with the blocked QR
algorithm \texttt{blockQR} (cf.~Figure \ref{fig:blockQR}), and then modifying it
by diagonalizing the diagonal blocks (as described in Section \ref{sec:diagonalizing}),
and then applying the high-accuracy pivoting scheme described in Section \ref{sec:power}.
The resulting algorithm \texttt{blockRRQR} is summarized in Figure \ref{fig:blockRRQR}.
The sparsity pattern of the matrix after each step of the algorithm is shown in Figure
\ref{fig:blockRRQRsparsity}.

\begin{figure}
\fbox{\begin{minipage}{150mm}
\begin{itemize}
\item Set $\mtx{Q} = \mtx{I}$ and $\mtx{P} = \mtx{I}$.
\item
\textbf{for} $i = 1,2,3,\dots,p-1$
\begin{itemize}
\item Partition the index vector $I = [I_{1}\ I_{2}\ I_{3}]$ so that $I_{2} = b(i-1) + (1:b)$ is the block processed in step $i$.
Then partition the matrix accordingly,
$$
\mtx{A} =
\left[\begin{array}{ccc}
\mtx{A}_{11} & \mtx{A}_{12} & \mtx{A}_{13} \\
\mtx{0}      & \mtx{A}_{22} & \mtx{A}_{23} \\
\mtx{0}      & \mtx{A}_{32} & \mtx{A}_{33}
\end{array}\right].
$$
(Observe that in the first step, $I_{1} = []$ and $\mtx{A}$ has only $2\times 2$ blocks.)
\item Set $\mtx{X} = \mtx{A}([I_{2}\ I_{3}],[I_{2}\ I_{3}])$. Then determine a pivoting matrix $\tilde{\mtx{P}}$ by processing
$(\mtx{X}^{*}\mtx{X})^{q}\mtx{X}^{*}$, as described in Section \ref{sec:power}.
Then update the last two block columns of $\mtx{A}$ accordingly,
$$
\left[\begin{array}{ccc}
\mtx{A}_{12}' & \mtx{A}_{13}' \\
\mtx{A}_{22}' & \mtx{A}_{23}' \\
\mtx{A}_{32}' & \mtx{A}_{33}'
\end{array}\right] =
\left[\begin{array}{ccc}
\mtx{A}_{12} & \mtx{A}_{13} \\
\mtx{A}_{22} & \mtx{A}_{23} \\
\mtx{A}_{32} & \mtx{A}_{33}
\end{array}\right]\,\tilde{\mtx{P}}.
$$
\item Execute a full SVD
$$
\begin{array}{ccccccccccc}
\vtwo{\mtx{A}_{22}'}{\mtx{A}_{23}'} &=& \tilde{\mtx{U}} & \vtwo{\mtx{D}_{22}}{\mtx{0}} & \tilde{\mtx{V}}^{*},\\
m'\times b && m'\times m' & m' \times b & b\times m'
\end{array}
$$
where $m' = m - (i-1)b$ is the length of the index vector $[I_{2}\ I_{3}]$. Observe that while $\tilde{\mtx{U}}$ is
large, it has internal structure that allows it to be applied efficiently, cf.~Section \ref{sec:diagonalizing}.
%(Specifically, it consists ``mostly'' of a product of $b$ Householder reflectors.)
\item Compute the new blocks $\mtx{A}_{23}''$ and $\mtx{A}_{33}''$ via
$\displaystyle
\vtwo{\mtx{A}_{23}''}{\mtx{A}_{23}''} = \tilde{\mtx{U}}^{*}
\vtwo{\mtx{A}_{23}'}{\mtx{A}_{23}'}.
$
\item Update the matrices $\mtx{A}$, $\mtx{Q}$, and $\mtx{P}$, via
\begin{align*}
\mtx{A}\leftarrow&\
\left[\begin{array}{ccc}
\mtx{A}_{11} & \mtx{A}_{12}'\tilde{\mtx{V}} & \mtx{A}_{13}'  \\
\mtx{0}      & \mtx{D}_{22}             & \mtx{A}_{23}'' \\
\mtx{0}      & \mtx{0}                  & \mtx{A}_{33}''
\end{array}\right],\\
\mtx{Q}\leftarrow&\ \mtx{Q}\left[\begin{array}{cc} \mtx{I}_{n_1} & \mtx{0} \\ \mtx{0} & \tilde{\mtx{U}}\end{array}\right],\\
\mtx{P}\leftarrow&\ \mtx{P}\left[\begin{array}{cc} \mtx{I}_{n_1} & \mtx{0} \\ \mtx{0} & \tilde{\mtx{P}}\end{array}\right]\,
\left[\begin{array}{ccc}
\mtx{I}_{n_1} & \mtx{0}         & \mtx{0} \\
\mtx{0}       & \tilde{\mtx{V}} & \mtx{0} \\
\mtx{0}       & \mtx{0}         & \mtx{I}_{n_3}
\end{array}\right].
\end{align*}
\end{itemize}
\textbf{end for}
\item At this point, all that remains is to process the lower right $b\times b$ block. Partition
$$
\mtx{A} = \mtwo{\mtx{A}_{11}}{\mtx{A}_{12}}{\mtx{0}}{\mtx{A}_{22}}
$$
so that $\mtx{A}_{22}$
is of size $b\times b$. Observe that $\mtx{A}_{11}$ is already upper triangular. Now compute the
(full, but small) SVD $\mtx{A}_{22} = \tilde{\mtx{U}}\mtx{D}_{22}\tilde{\mtx{V}}^{*}$. Then update
$$
\mtx{A}\leftarrow\mtwo{\mtx{A}_{11}}{\mtx{A}_{12}}{\mtx{0}}{\mtx{D}_{22}},\qquad
\mtx{Q} = \mtx{Q}\mtwo{\mtx{I}_{n_1}}{\mtx{0}}{\mtx{0}}{\tilde{\mtx{U}}},\qquad
\mtx{P} = \mtx{P}\mtwo{\mtx{I}_{n_1}}{\mtx{0}}{\mtx{0}}{\tilde{\mtx{V}}}.
$$
\end{itemize}
\end{minipage}}
\caption{The algorithm \texttt{blockRRQR}. Given an input matrix $\mtx{A}$ of size $m\times n$, with $m\geq n$,
and a block size $b$, the algorithm produces a factorization $\mtx{A}\mtx{P} = \mtx{Q}\mtx{R}$ with $\mtx{R}$
upper triangular, and $\mtx{P}$ and $\mtx{Q}$ orthonormal, cf.~(\ref{eq:AP=QR}). For simplicity, we assume
that $n = bp$ for some integer $p$. The matrix $\mtx{R}$ overwrites $\mtx{A}$.}
\label{fig:blockRRQR}
\end{figure}

\begin{figure}
\begin{tabular}{ccccc}
\includegraphics[height=32mm]{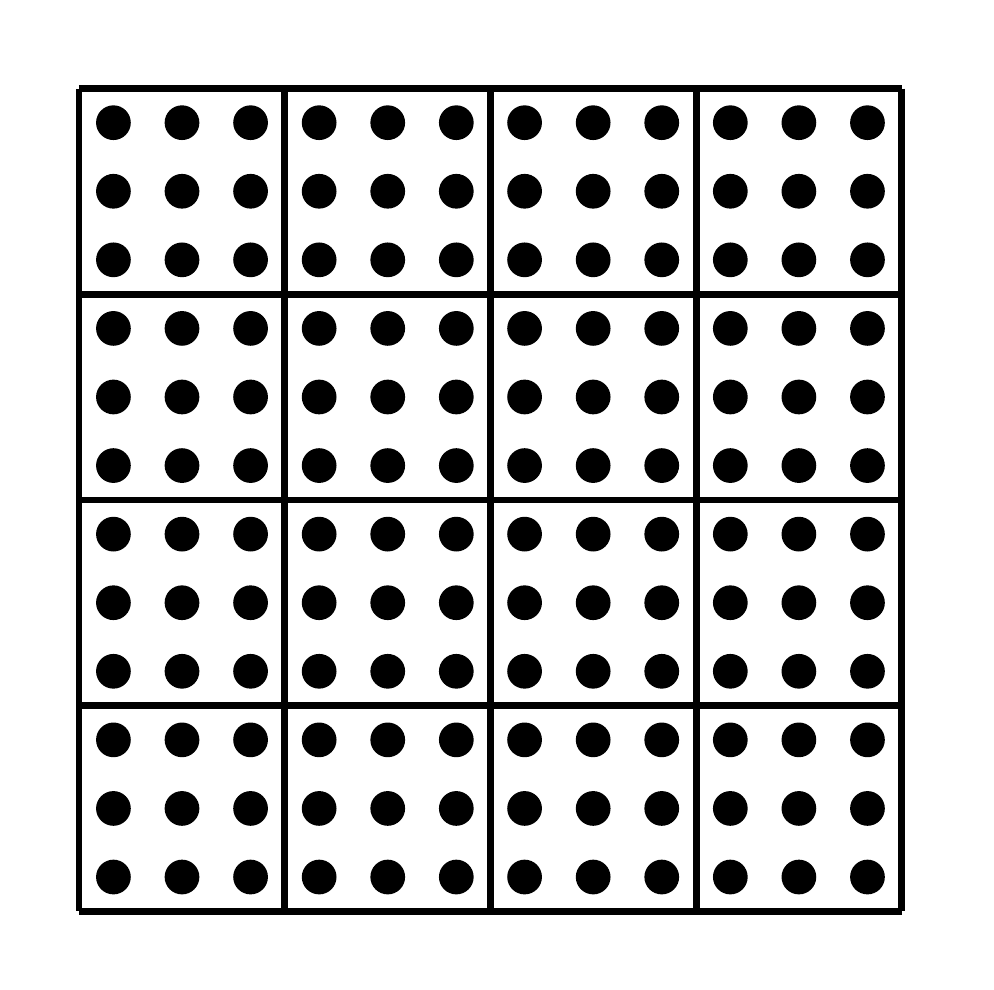} &
\includegraphics[height=32mm]{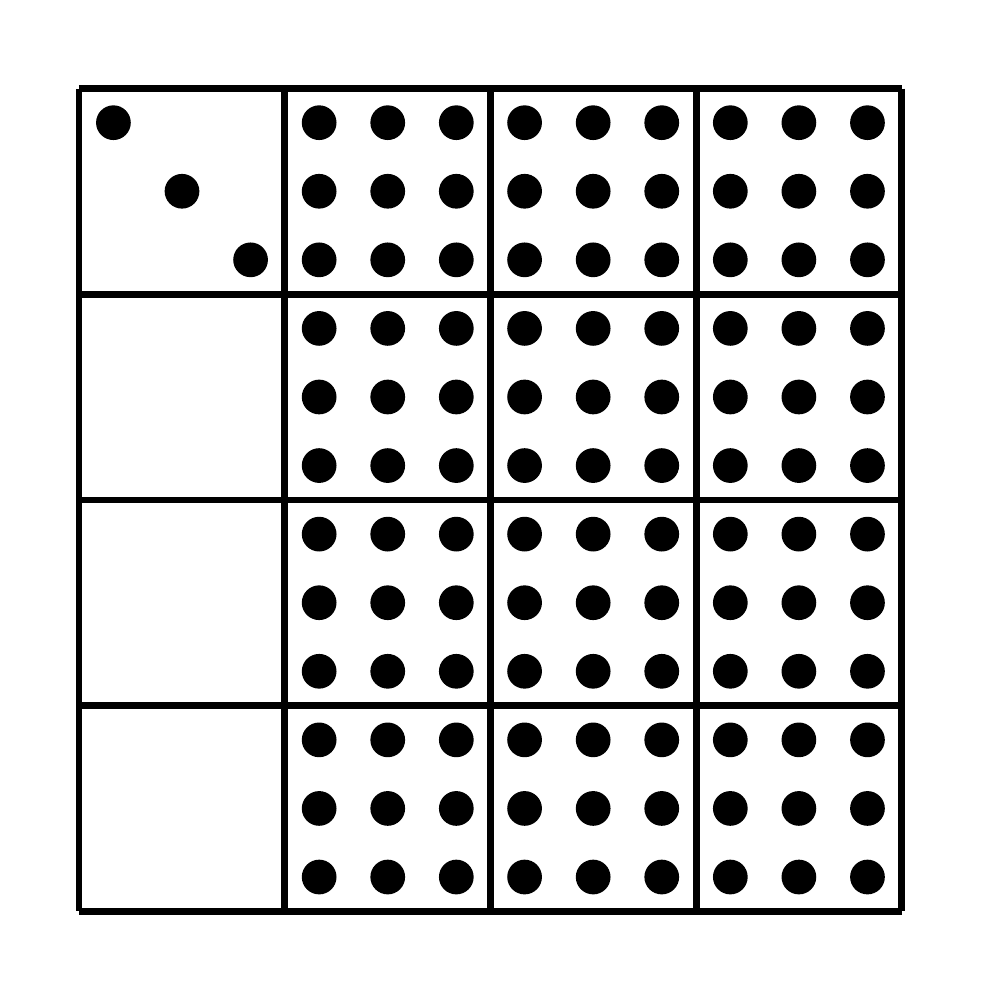} &
\includegraphics[height=32mm]{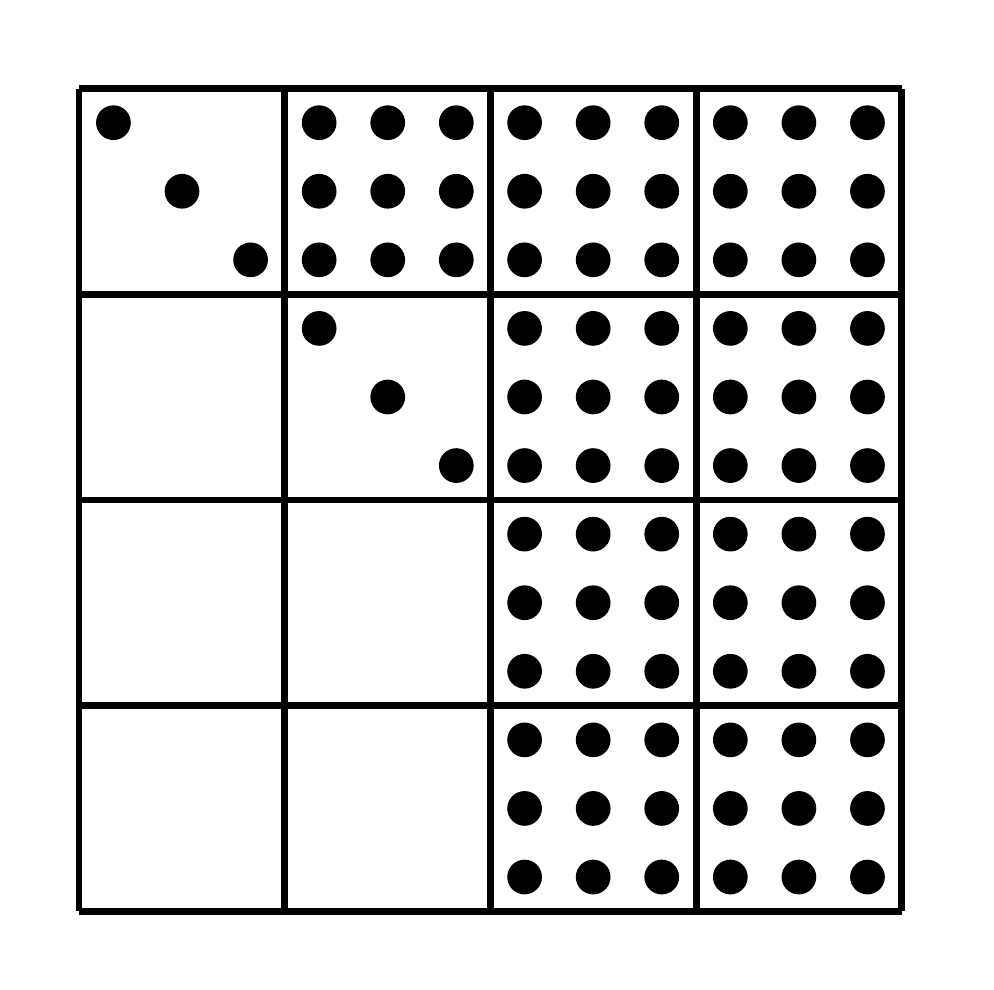} &
\includegraphics[height=32mm]{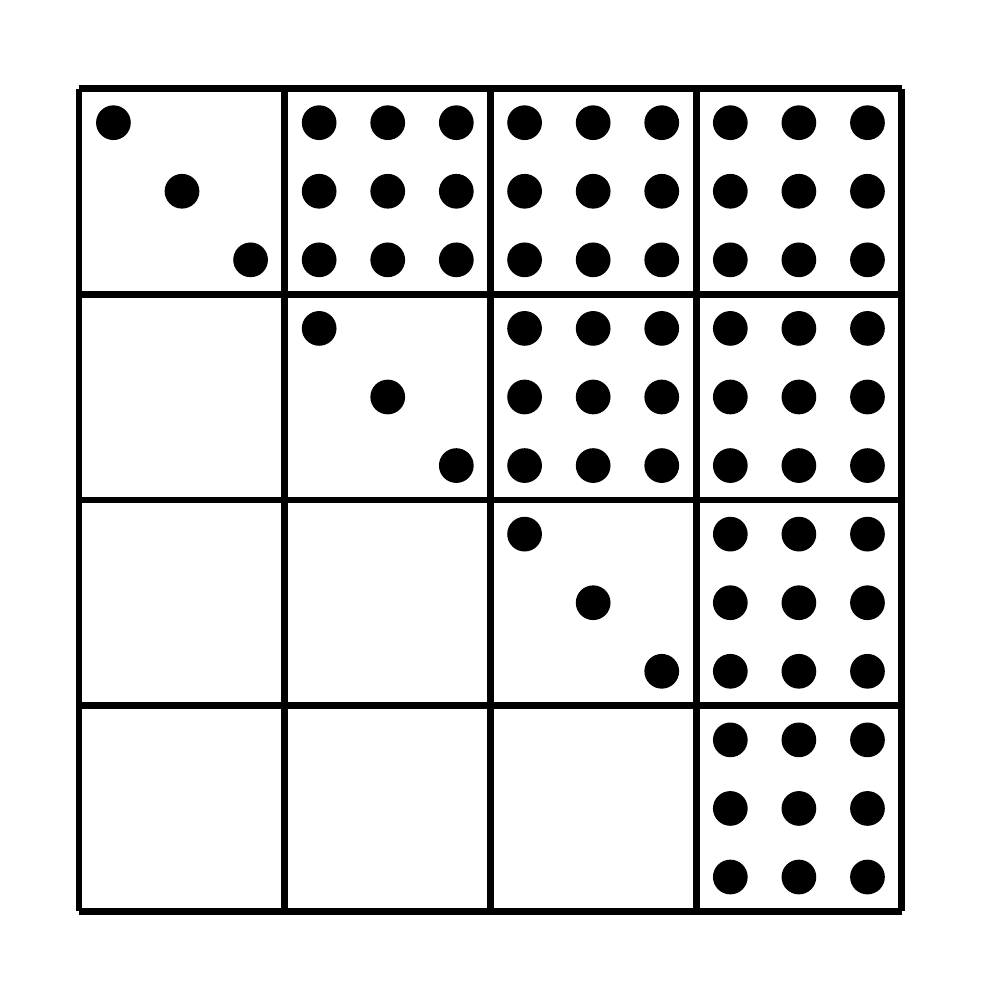} &
\includegraphics[height=32mm]{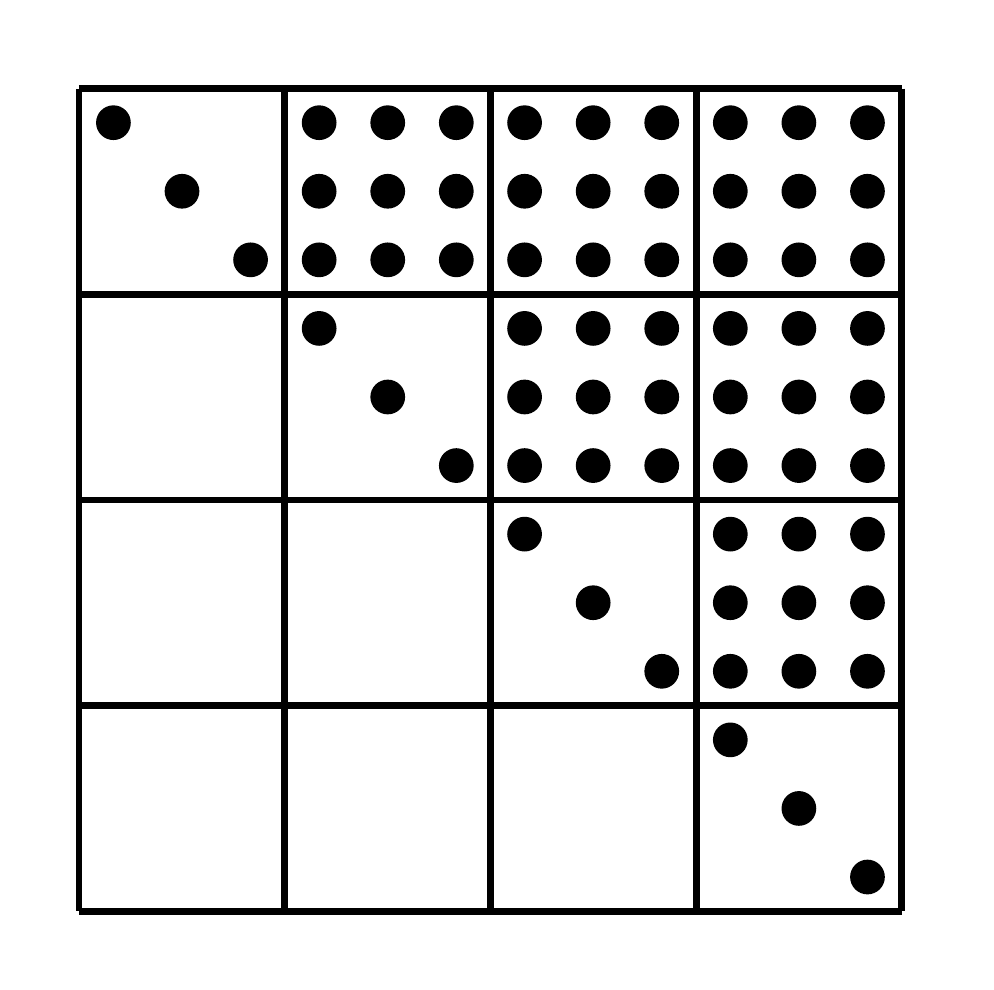} \\
$\mtx{A}_{0}=\mtx{A}$ &
$\mtx{A}_{1}$ &
$\mtx{A}_{2}$ &
$\mtx{A}_{3}$ &
$\mtx{A}_{4}=\mtx{R}$
\end{tabular}
\caption{Example of how a $12\times 12$ matrix $\mtx{A}$ is driven to upper triangular form in
the randomized RRQR factorization. Each figure shows the sparsity pattern of matrix $\mtx{A}_{j}$,
resulting after $j$ steps of the algorithm \texttt{blockRRQR}, cf.~Figure \ref{fig:blockRRQR}.}
\label{fig:blockRRQRsparsity}
\end{figure}

\clearpage

\section{Numerical experiments}

In this section, we test three different methods for computing a full QR factorization of
a given matrix $\mtx{A}$:
\begin{itemize}
\item[Method 1] The Algorithm \texttt{blockQR} as shown in Figure \ref{fig:blockQR}, with the ``pivoting matrix''
chosen as a permutation matrix as described in Section \ref{sec:pivotperm}. No over-sampling is used.
\item[Method 2] The Algorithm \texttt{blockQR} as shown in Figure \ref{fig:blockQR}, with the ``pivoting matrix''
chosen as a product of Householder reflectors as described in Section \ref{sec:pivotON}. No over-sampling is used.
\item[Method 3] The Algorithm \texttt{blockRRQR} as shown in Figure \ref{fig:blockRRQR}, with the over-sampling
parameter $p$ set to be half the block size $b$.
\end{itemize}
We did not include over-sampling when running \texttt{blockQR} since numerical experiments
indicated that there was essentially no benefit to doing so.
For simplicity, $\mtx{A}$ is in every experiment a real matrix of size $n\times n$.

At the time of writing, we have not yet implemented an optimized version of the algorithm, so we
cannot present timing comparisons. The purpose of the numerical experiments is to demonstrate the
very high accuracy of the proposed techniques.

\subsection{A matrix with rapidly decaying singular values}
\label{sec:numfast}
In our first experiments, we apply the various factorization algorithms to the matrix
$$
\mtx{A}^{\rm fast} = \mtx{U}\,\mtx{D}\,\mtx{V}^{*},
$$
with $\mtx{U}$ and $\mtx{V}$ unitary matrices drawn from a uniform distribution. The
matrix $\mtx{D}$ is diagonal, with diagonal entries
$$
\mtx{D}(i,i) = \bigl(10^{-5}\bigr)^{(i-1)/(n-1)}.
$$
In other words, the singular values of $\mtx{A}^{\rm fast}$ decay exponentially, from $\sigma_{1} = 1$ to $\sigma_{n} = 10^{-5}$.
In the experiments shown, we set $n=300$, and use the block size $b=50$.

\begin{figure}
\includegraphics[width=0.95\textwidth]{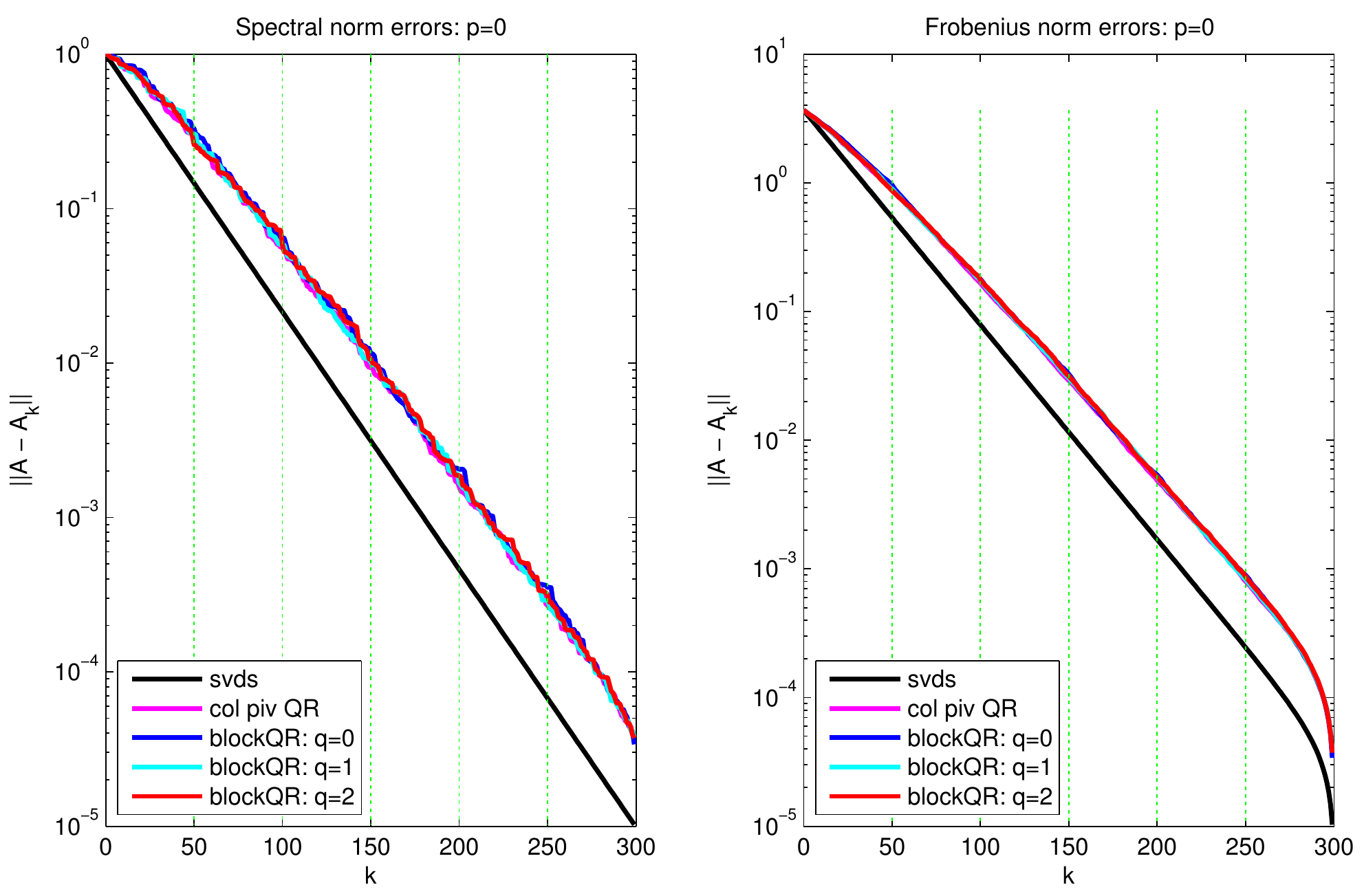}
\caption{Approximation errors from Method 1 (\texttt{blockQR} with a permutation matrix) for
$\mtx{A} = \mtx{A}^{\rm fast}$, cf.~Section \ref{sec:numfast}.
Also included are the errors resulting from a truncated SVD (which are theoretically optimal), and
a column pivoted QR factorization.}
\label{fig:err_method1_matrix1}
\end{figure}

\begin{figure}
\includegraphics[width=110mm]{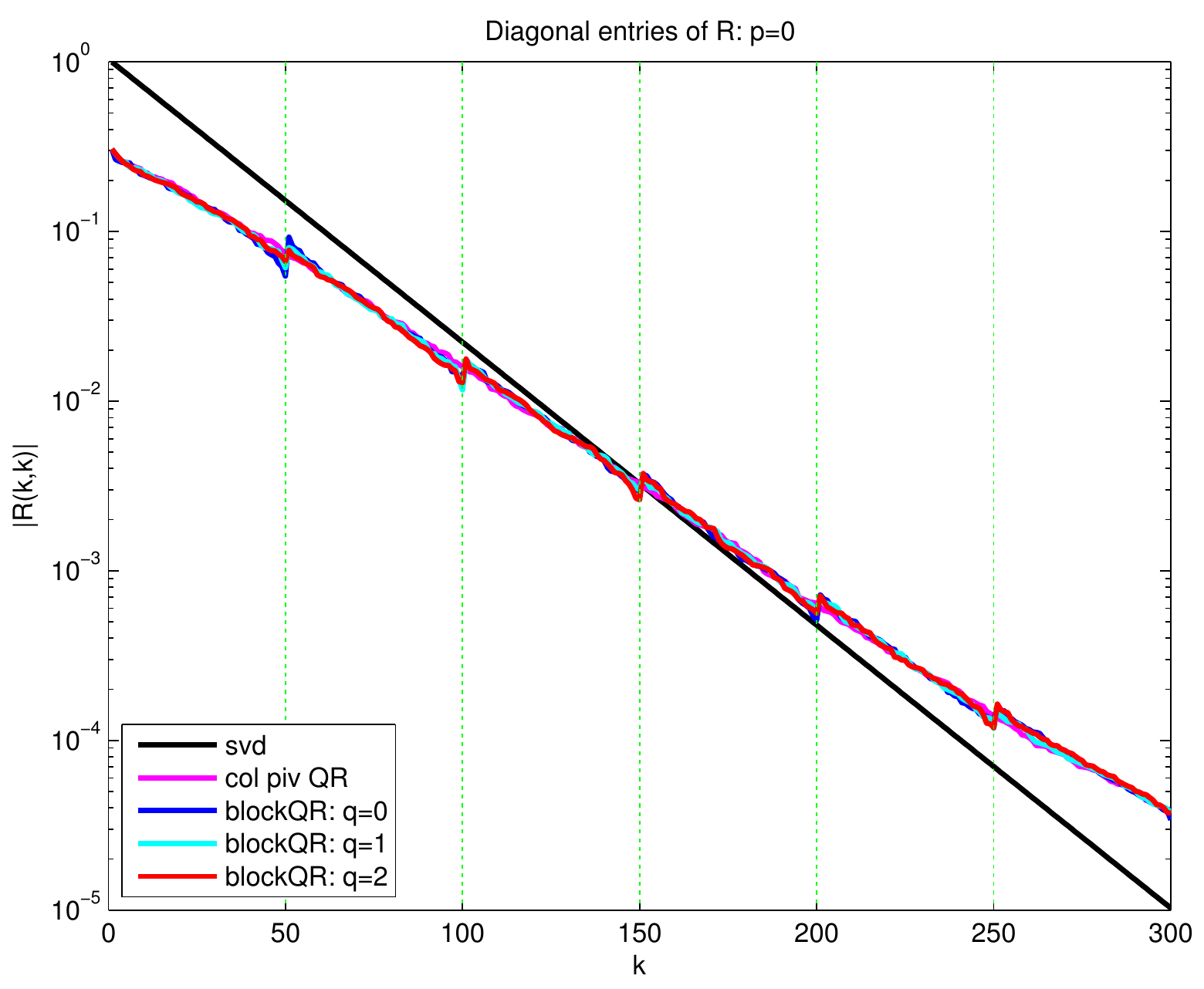}
\caption{Comparison of $|\mtx{R}(k,k)|$ for the matrix $\mtx{R}$ resulting from Method 1,
and the singular values of $\mtx{A}^{\rm fast}$. We would ideally like these values to be close.}
\label{fig:diag_method1_matrix1}
\end{figure}

We first study ``Method 1'' (\texttt{blockQR} with a pivoting matrix).
Figure \ref{fig:err_method1_matrix1} shows the error $e_{k}$ obtained when truncating the factorization,
$$
e_{k} = \|\mtx{A}^{\rm fast} - \mtx{A}_{k}^{\rm fast}\|,
\qquad\mbox{where}\qquad
\mtx{A}_{k}^{\rm fast} = \mtx{Q}(:,1:k)\,\mtx{R}(1:k,:)\,\mtx{P}^{*}.
$$
The figure also shows the corresponding errors when $\mtx{A}_{k}^{\rm fast}$ is the truncated SVD
and the truncated QR factorization obtained from classical column pivoting. We make three
observations:
\begin{itemize}
\item The errors resulting from Method 1 are very close to the errors obtained by the column pivoted QR factorization.
\item The errors from every truncated QR factorization involving a permutation matrix that we tried are substantially sub-optimal,
as compared to the truncated SVD.
\item Using the ``power method'' described in Section \ref{sec:power} leads to almost no improvement in this case.
\end{itemize}
Figure \ref{fig:diag_method1_matrix1} shows how the diagonal entries of $\{|\mtx{R}(k,k)|\}_{k=1}^{n}$ compare to the
singular values $\{\mtx{D}(k,k)\}_{k=1}^{n}$. The figure shows that Method 1 performs no better (and no worse) than
classical column pivoted QR.

\begin{figure}
\includegraphics[width=0.95\textwidth]{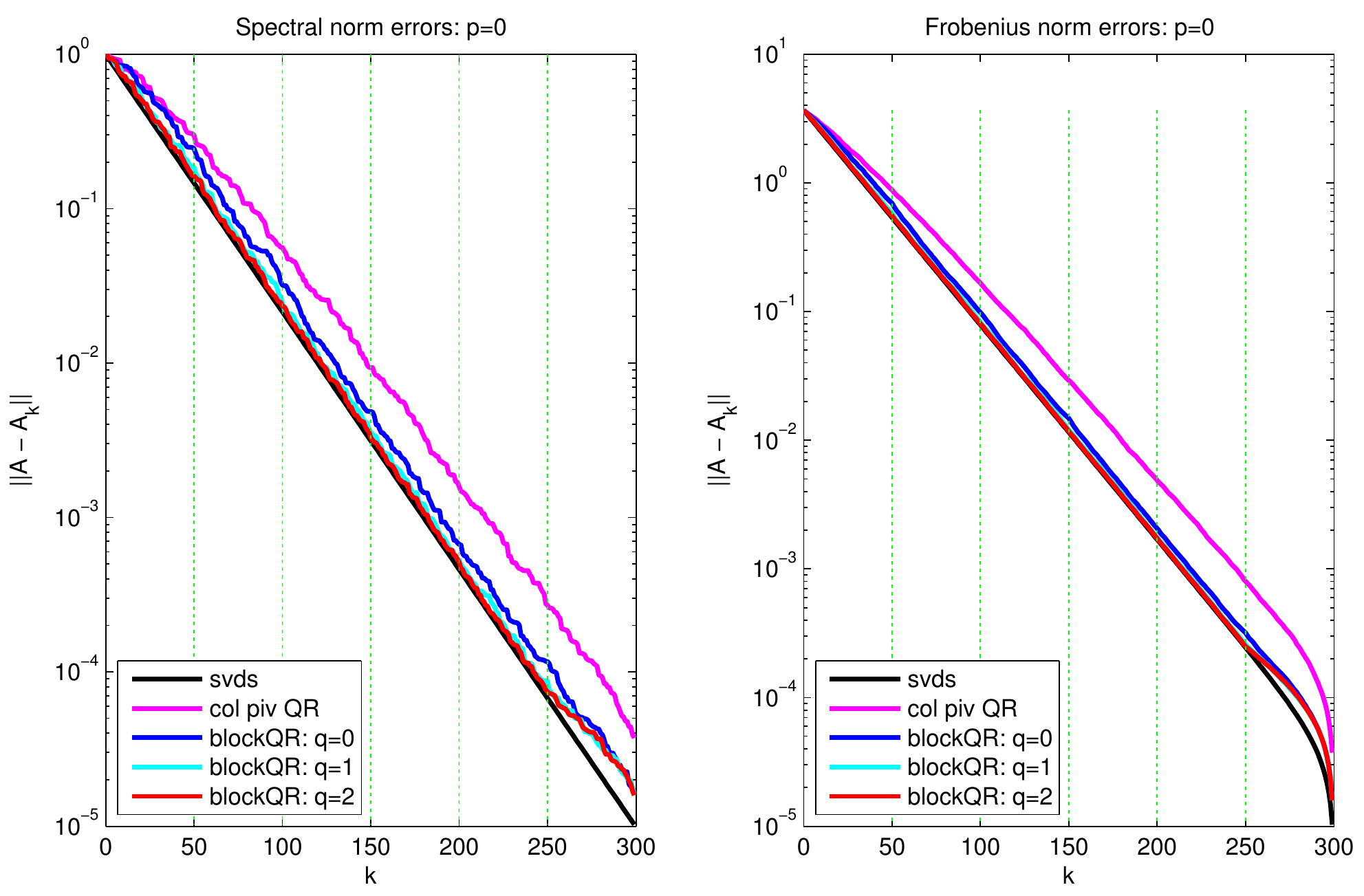}
\caption{Approximation errors from Method 2 (\texttt{blockQR} with a Householder pivoting matrix) for
$\mtx{A} = \mtx{A}^{\rm fast}$, cf.~Section \ref{sec:numfast}.
Also included are the errors resulting from a truncated SVD (which are theoretically optimal), and
a column pivoted QR factorization.}
\label{fig:err_method2_matrix1}
\end{figure}

\begin{figure}
\includegraphics[width=110mm]{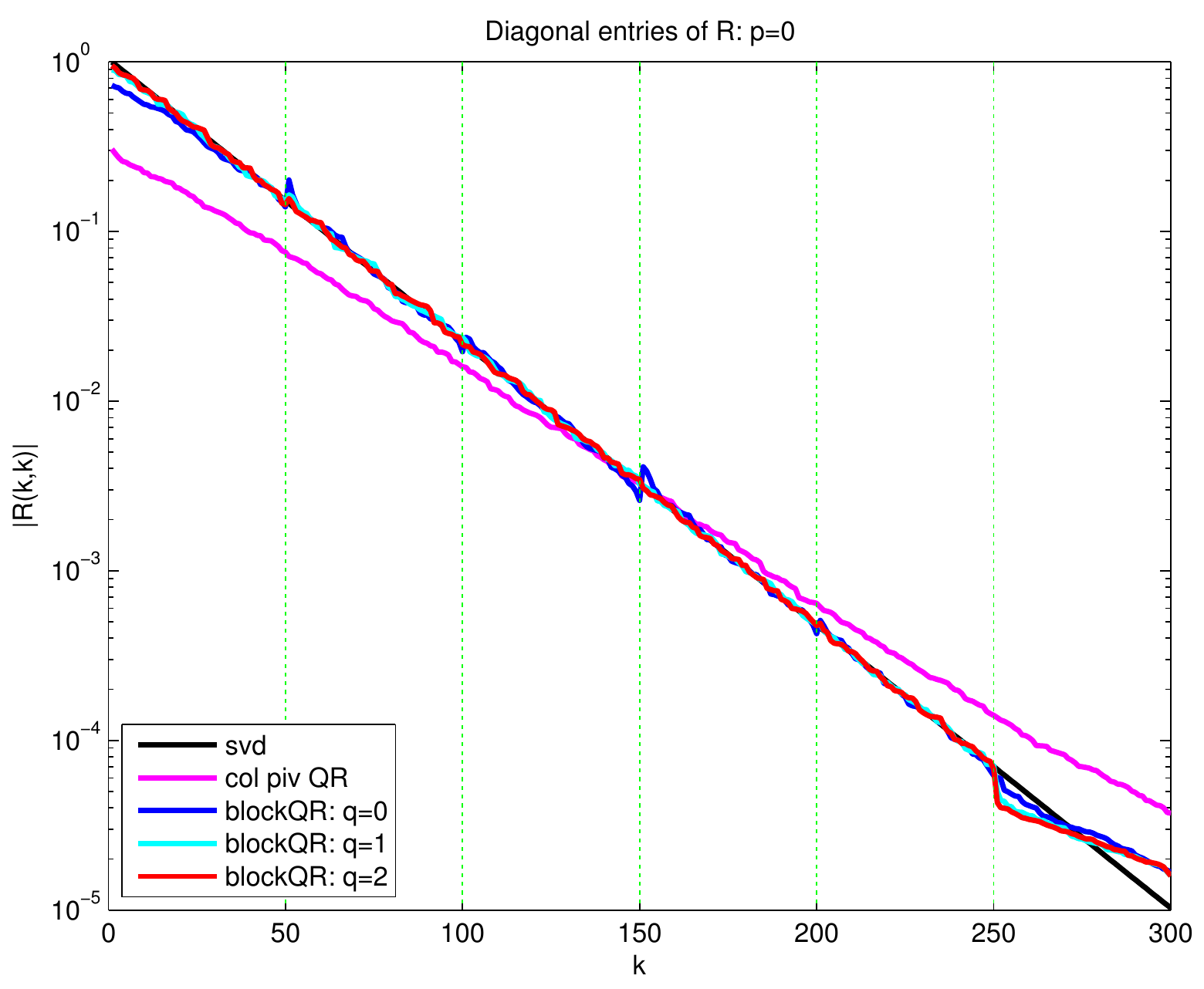}
\caption{Comparison of $|\mtx{R}(k,k)|$ for the matrix $\mtx{R}$ resulting from Method 2,
and the singular values of $\mtx{A}^{\rm fast}$.}
\label{fig:diag_method2_matrix1}
\end{figure}

We next use ``Method 2'' (\texttt{blockQR} with a Householder pivoting matrix) to compute the factorization
(\ref{eq:AP=QR}), with the results shown in Figures \ref{fig:err_method2_matrix1} and \ref{fig:diag_method2_matrix1}.
The approximation error is now much better even without using the power method, and once the power method is
employed, results improve very rapidly.

Finally, we test ``Method 3'' (\texttt{blockRRQR} with over-sampling parameter $p=25$) to compute the factorization
(\ref{eq:AP=QR}), with the results shown in Figures \ref{fig:err_method2_matrix1} and \ref{fig:diag_method2_matrix1}.
The approximation error is now much better even without using the power method, and once the power method is
employed, results improve very rapidly.

\begin{figure}
\includegraphics[width=0.95\textwidth]{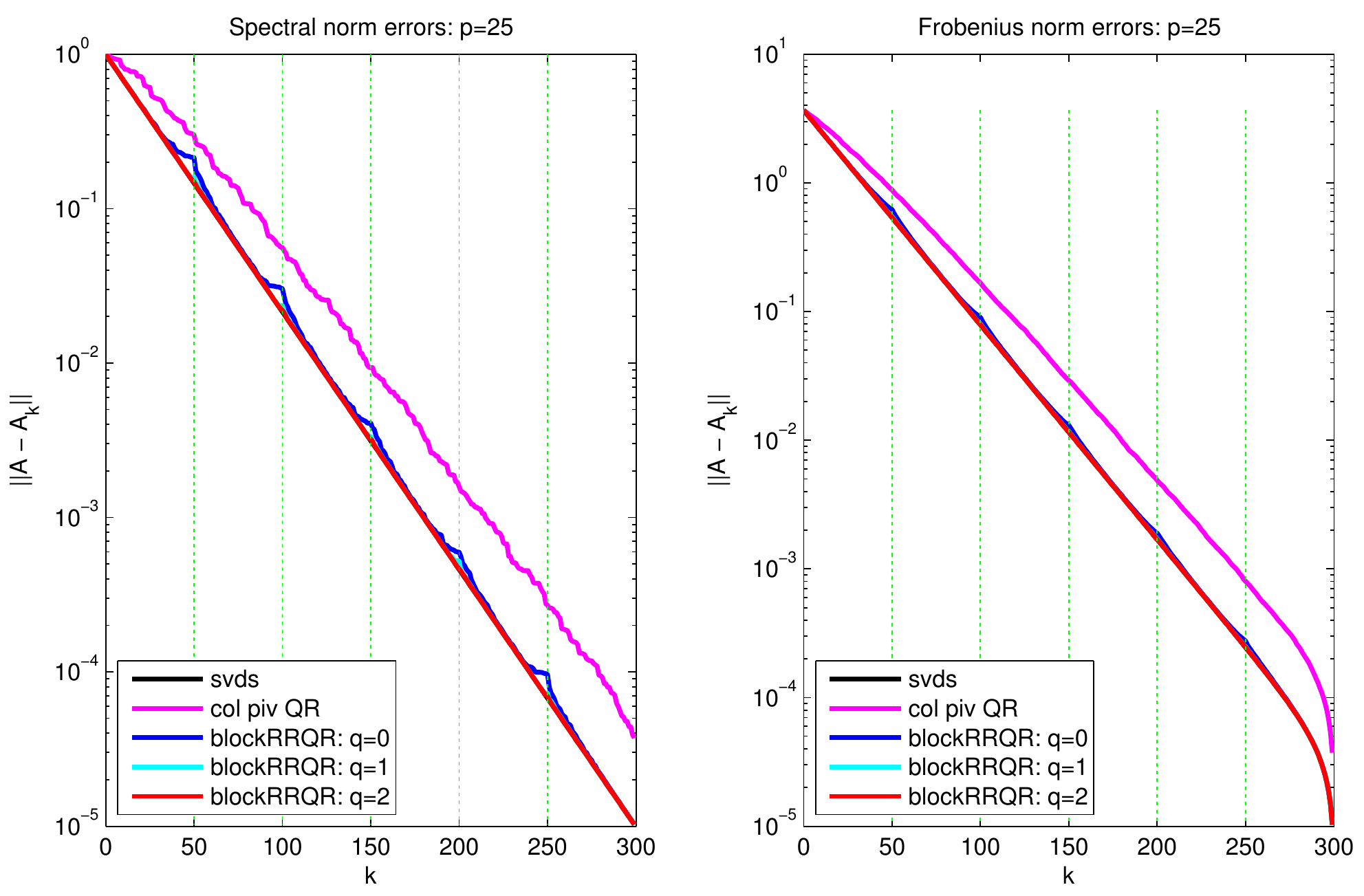}
\caption{Approximation errors from Method 3 (\texttt{blockRRQR} with $p=25$) for
$\mtx{A} = \mtx{A}^{\rm fast}$, cf.~Section \ref{sec:numfast}.
Also included are the errors resulting from a truncated SVD (which are theoretically optimal), and
a column pivoted QR factorization.}
\label{fig:err_method3_matrix1}
\end{figure}

\begin{figure}
\includegraphics[width=110mm]{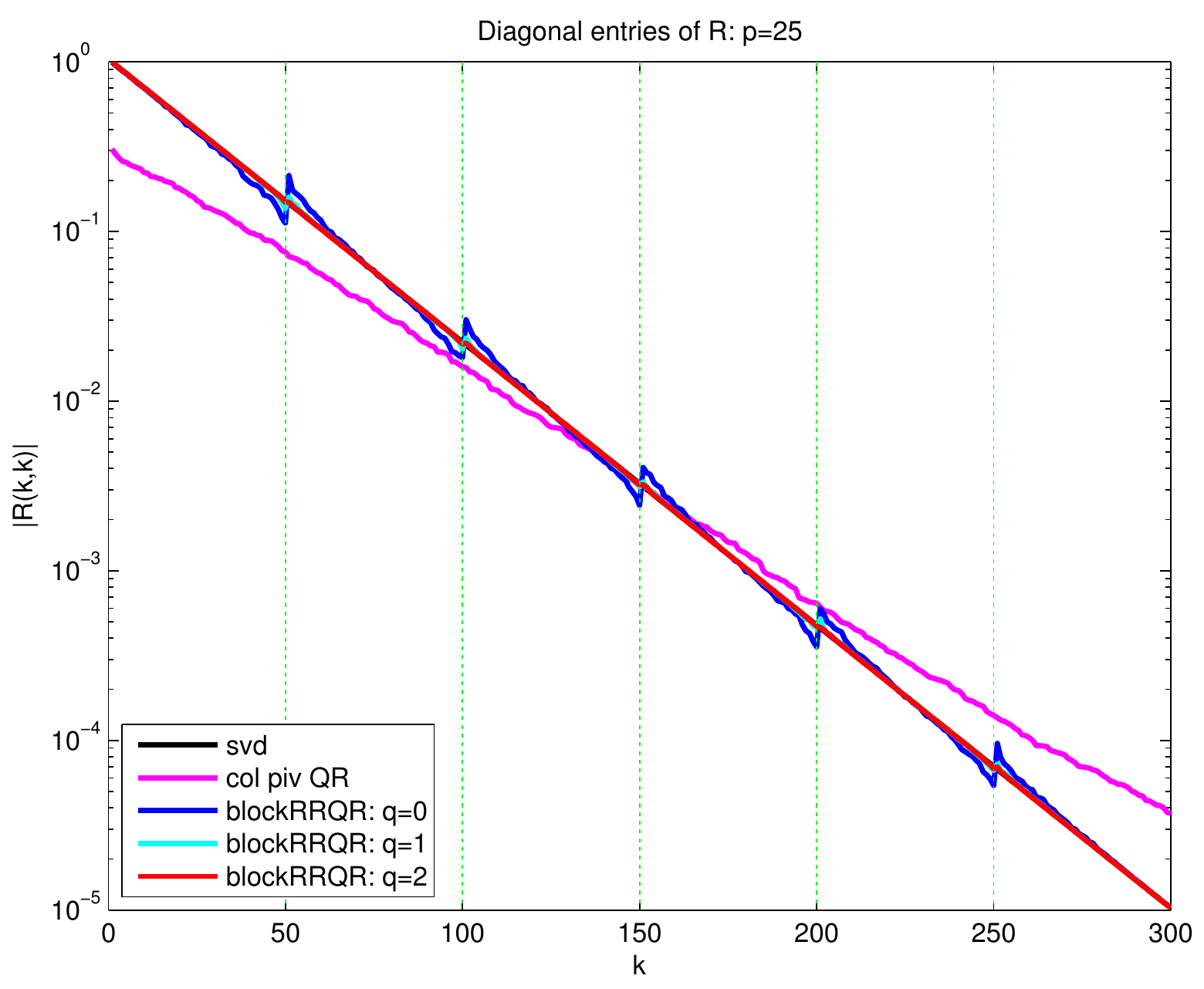}
\caption{Comparison of $|\mtx{R}(k,k)|$ for the matrix $\mtx{R}$ resulting from Method 3,
and the singular values of $\mtx{A}^{\rm fast}$.}
\label{fig:diag_method3_matrix1}
\end{figure}

In all experiments shown, we see that the randomized methods are particularly good at
minimizing errors in the Frobenius norm.

\clearpage
\subsection{A random Gaussian matrix}
\label{sec:numgauss}
\begin{figure}
\includegraphics[width=0.95\textwidth]{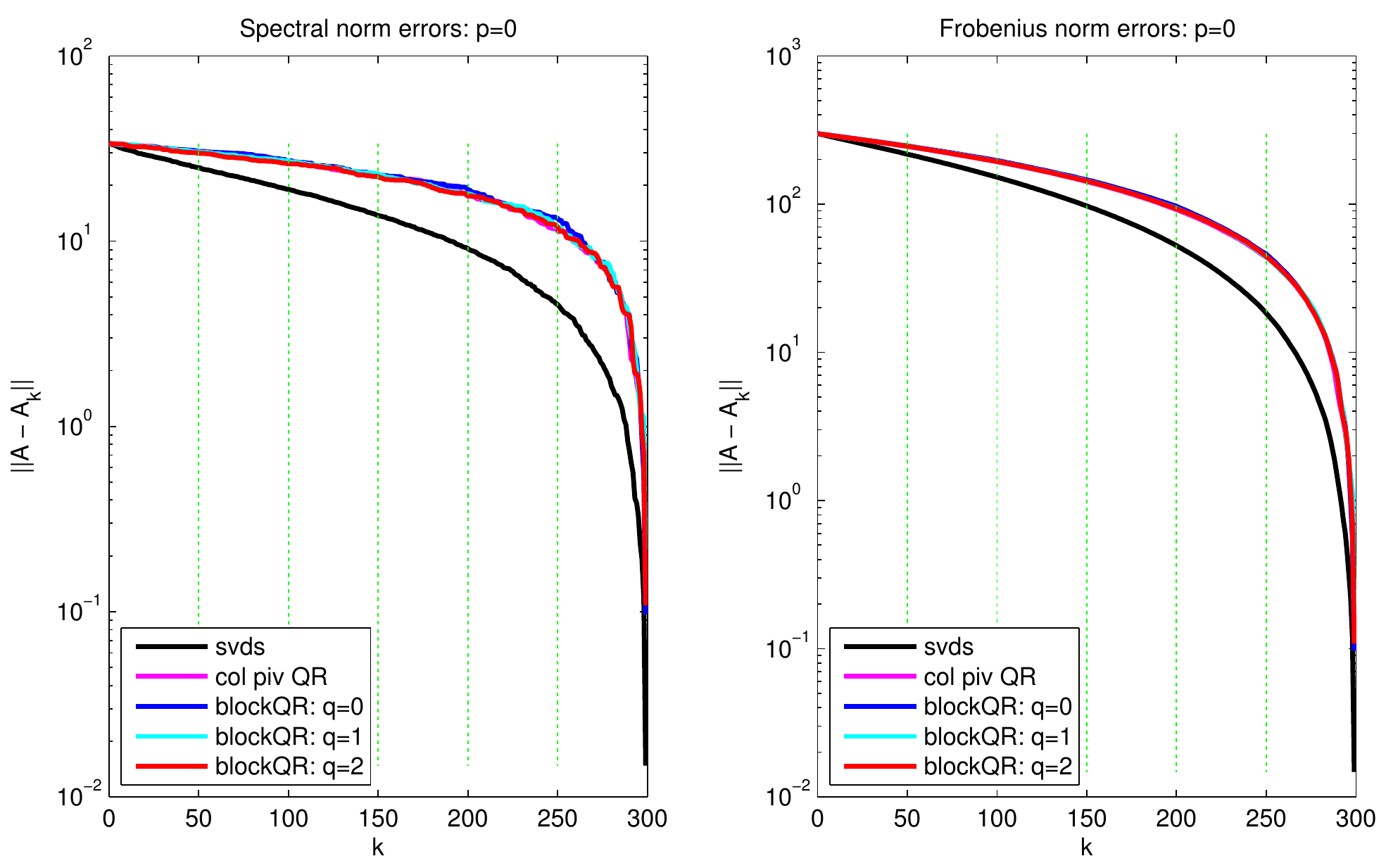}
\caption{Approximation errors from Method 1 (\texttt{blockQR} with a permutation matrix) for
$\mtx{A} = \mtx{A}^{\rm gauss}$, cf.~Section \ref{sec:numgauss}.
Also included are the errors resulting from a truncated SVD (which are theoretically optimal), and
a column pivoted QR factorization.}
\label{fig:err_method1_matrix2}
\end{figure}

\begin{figure}
\includegraphics[width=110mm]{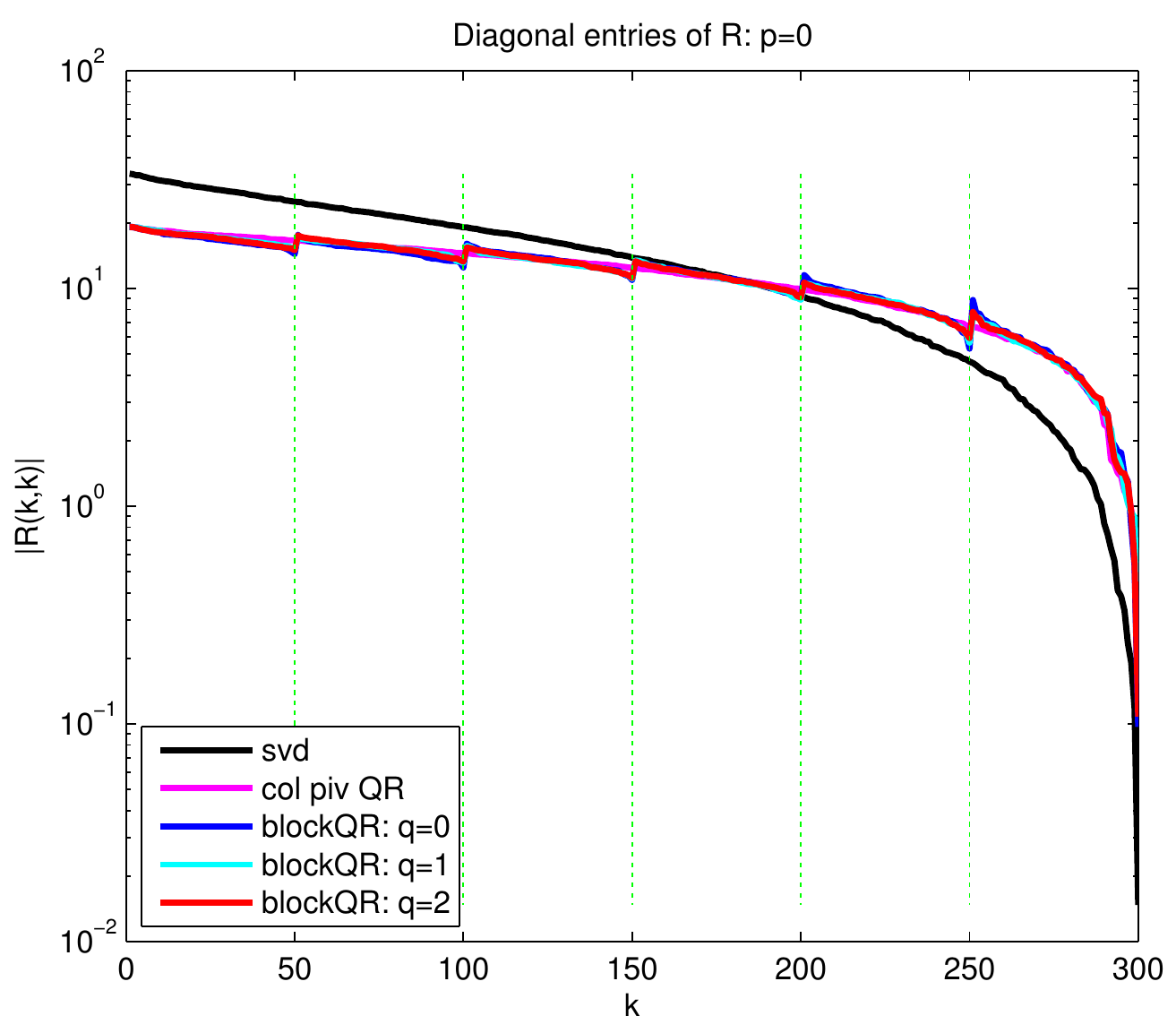}
\caption{Comparison of $|\mtx{R}(k,k)|$ for the matrix $\mtx{R}$ resulting from Method 1,
and the singular values of $\mtx{A}^{\rm gauss}$. We would ideally like these values to be close.}
\label{fig:diag_method1_matrix2}
\end{figure}

\begin{figure}
\includegraphics[width=0.95\textwidth]{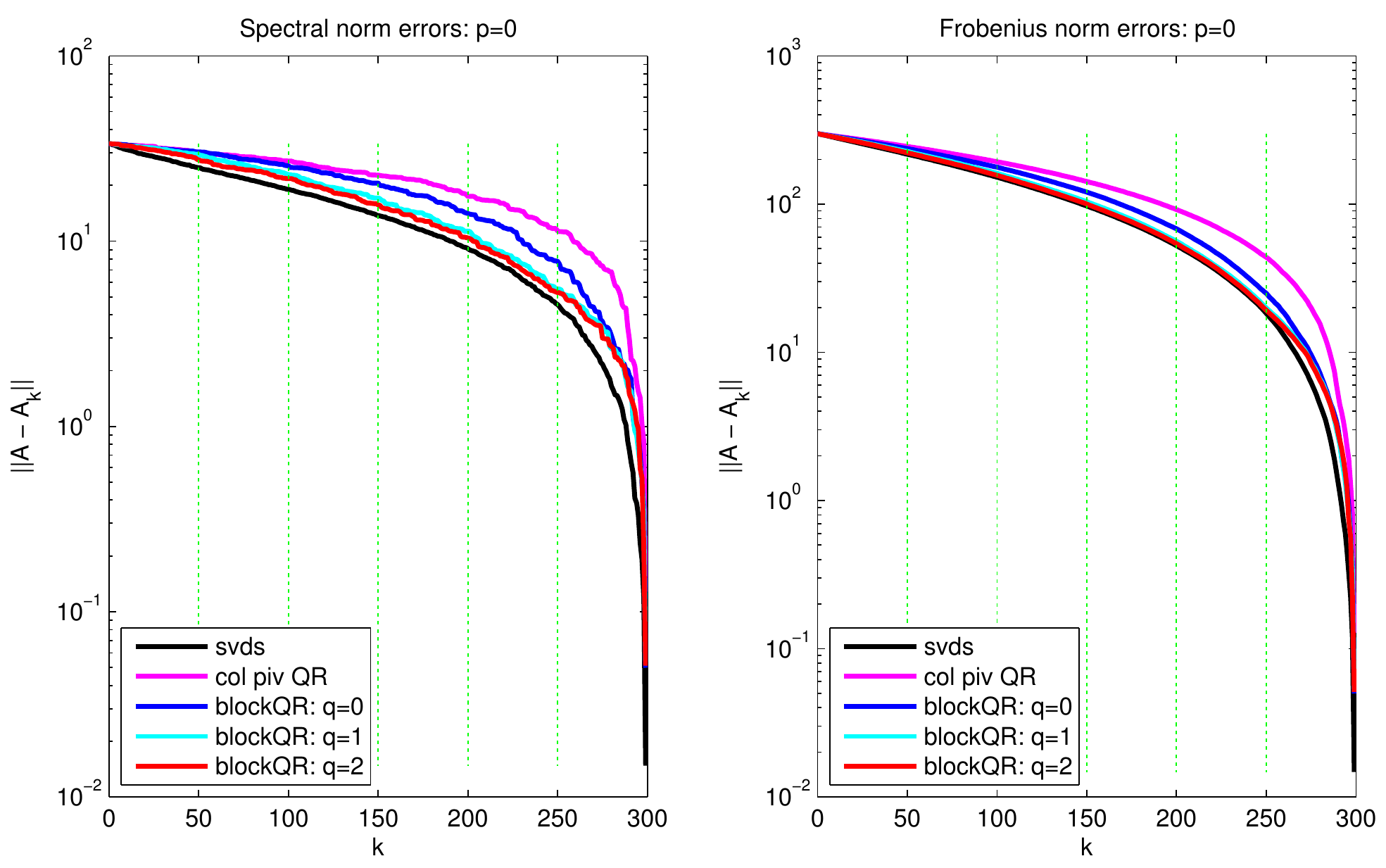}
\caption{Approximation errors from Method 2 (\texttt{blockQR} with a Householder pivoting matrix) for
$\mtx{A} = \mtx{A}^{\rm gauss}$, cf.~Section \ref{sec:numgauss}.
Also included are the errors resulting from a truncated SVD (which are theoretically optimal), and
a column pivoted QR factorization.}
\label{fig:err_method2_matrix2}
\end{figure}

\begin{figure}
\includegraphics[width=110mm]{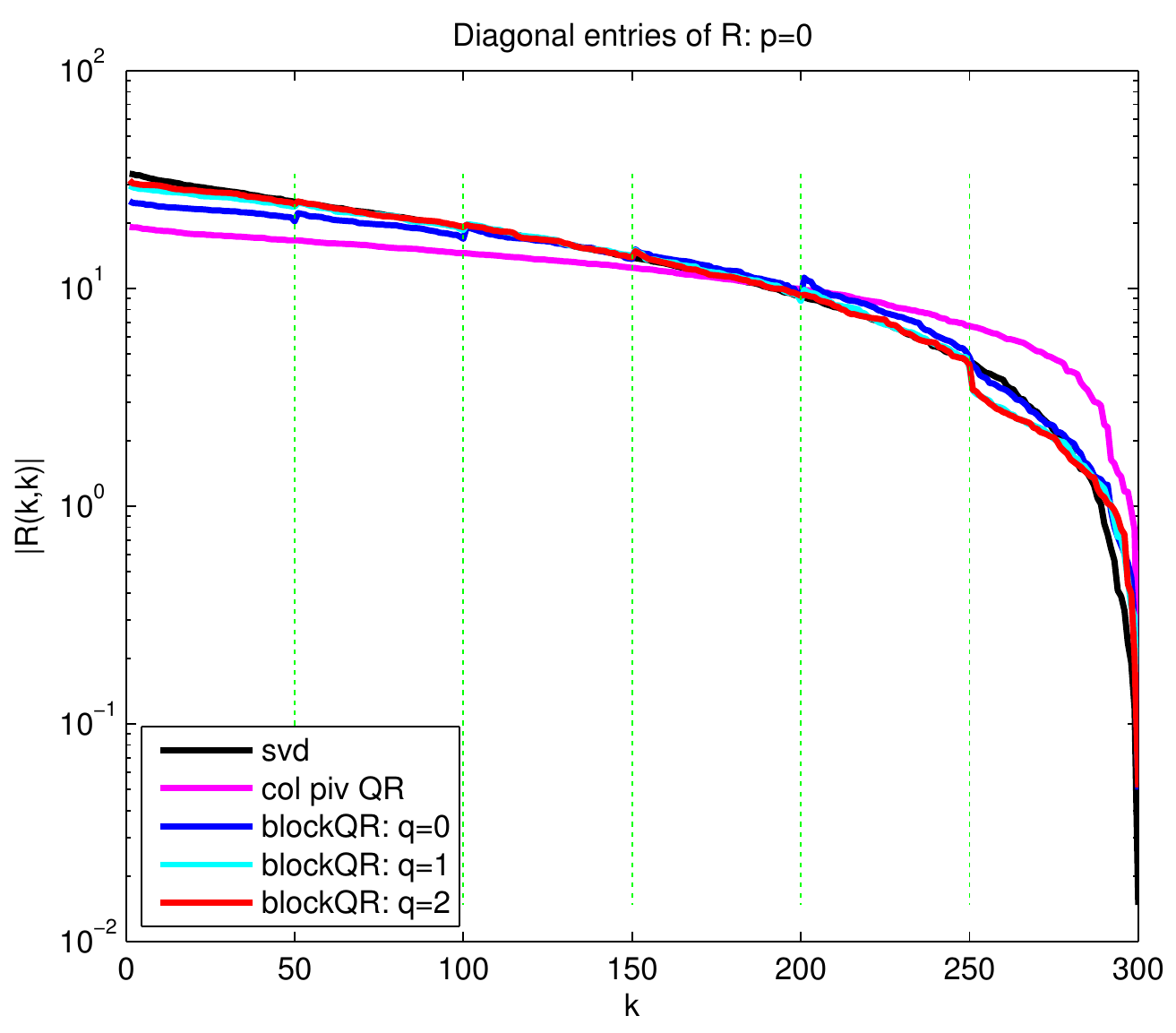}
\caption{Comparison of $|\mtx{R}(k,k)|$ for the matrix $\mtx{R}$ resulting from Method 2,
and the singular values of $\mtx{A}^{\rm gauss}$.}
\label{fig:diag_method2_matrix2}
\end{figure}

\begin{figure}
\includegraphics[width=0.95\textwidth]{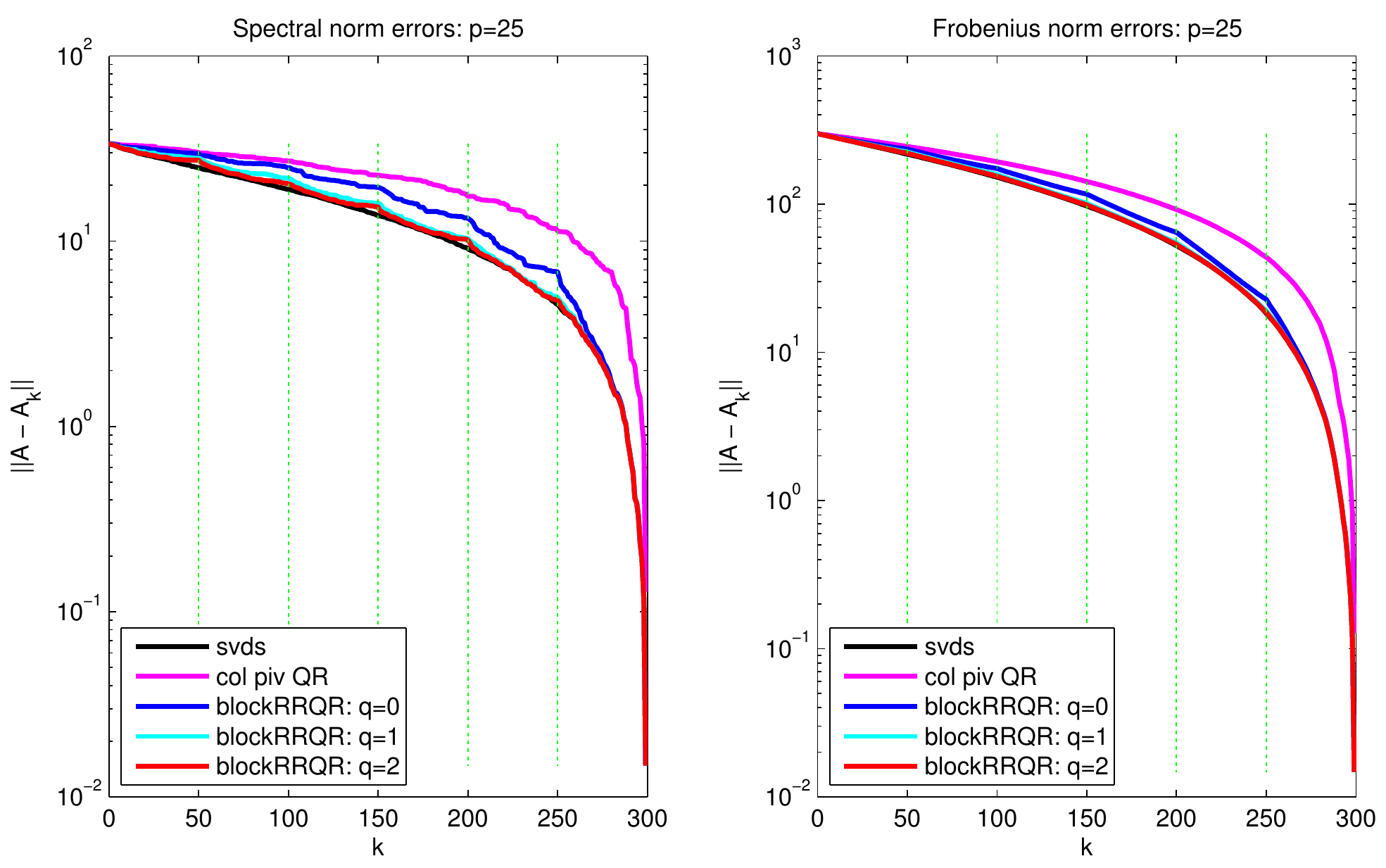}
\caption{Approximation errors from Method 3 (\texttt{blockRRQR} with $p=25$) for
$\mtx{A} = \mtx{A}^{\rm gauss}$, cf.~Section \ref{sec:numgauss}.
Also included are the errors resulting from a truncated SVD (which are theoretically optimal), and
a column pivoted QR factorization.}
\label{fig:err_method3_matrix2}
\end{figure}

\begin{figure}
\includegraphics[width=110mm]{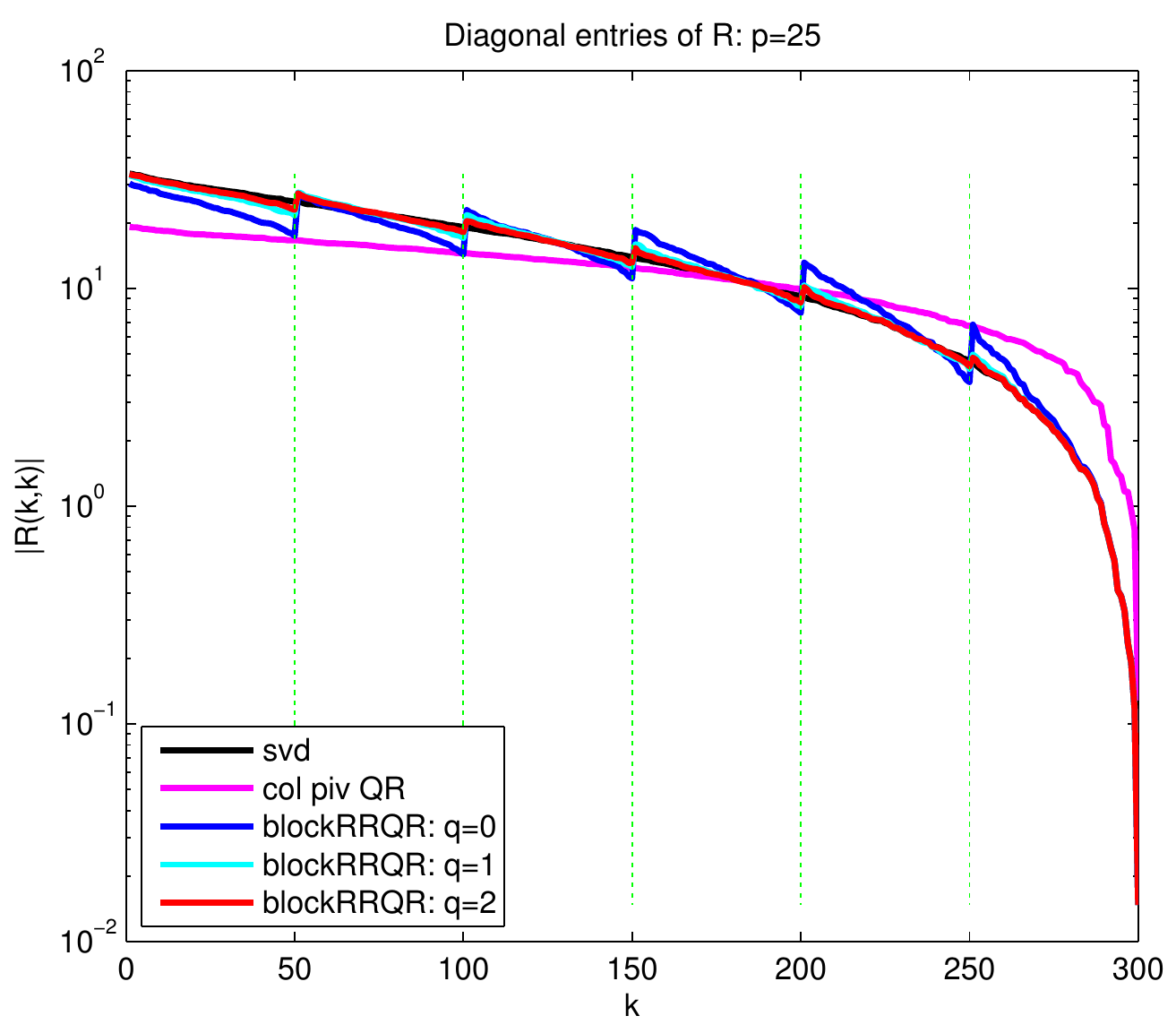}
\caption{Comparison of $|\mtx{R}(k,k)|$ for the matrix $\mtx{R}$ resulting from Method 3,
and the singular values of $\mtx{A}^{\rm gauss}$.}
\label{fig:diag_method3_matrix2}
\end{figure}

We next repeat all experiments conducted in Section \ref{sec:numfast}, but now for a
square matrix $\mtx{A}^{\rm gauss}$ with i.i.d.~normalized Gaussian entries. The singular values of this matrix
initially decay slowly, but then plummet at the end. The matrix is again of size $300\times 300$,
and we used a block size of $b=50$. The results are shown in Figures \ref{fig:err_method1_matrix2} -- \ref{fig:diag_method3_matrix2}.

We see that the results for $\mtx{A}^{\rm guass}$, whose singular values decay slowly, are quite
similar to those for $\mtx{A}^{\rm fast}$, whose singular values decay rapidly. The performance
of Method 1 is again very similar to that of classical column pivoted QR, with very little benefit
seen from using the power method. Once we allow the permutation matrix to be Householder reflectors,
the errors improve greatly, in particular once the power method is employed. For this example, it
is worth noting that for Method 3 with the power method engaged with $q=2$, the diagonal entries of
$\mtx{R}$ are \textit{very} close to the true singular values, even at the very end of the spectrum.

\clearpage

\clearpage
\subsection{A matrix with S-shaped decay in its singular values}
\label{sec:numS}
\begin{figure}
\includegraphics[width=0.95\textwidth]{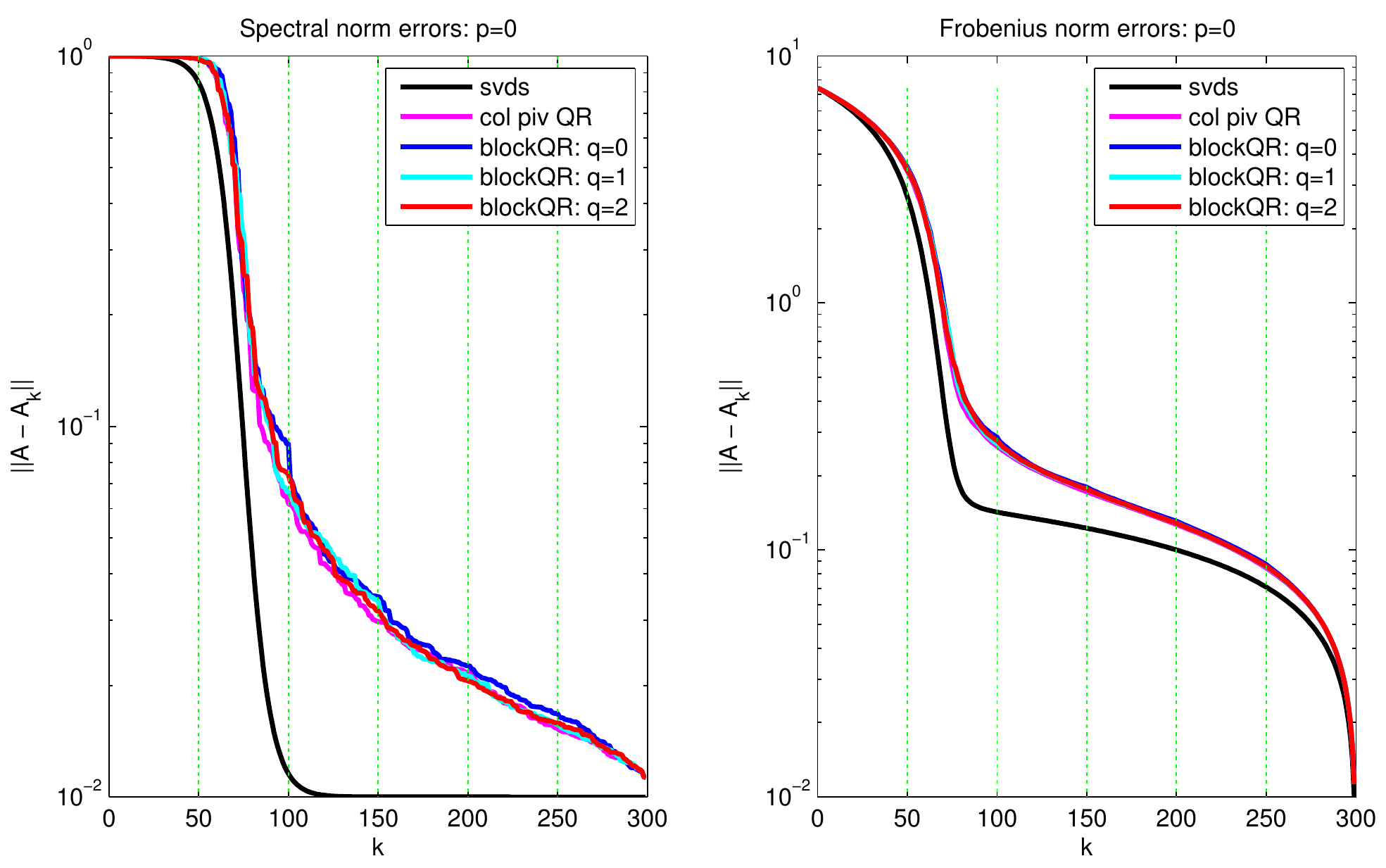}
\caption{Approximation errors from Method 1 (\texttt{blockQR} with a permutation matrix) for
$\mtx{A} = \mtx{A}^{\rm S}$, cf.~Section \ref{sec:numS}.
Also included are the errors resulting from a truncated SVD (which are theoretically optimal), and
a column pivoted QR factorization.}
\label{fig:err_method1_matrix3}
\end{figure}

\begin{figure}
\includegraphics[width=110mm]{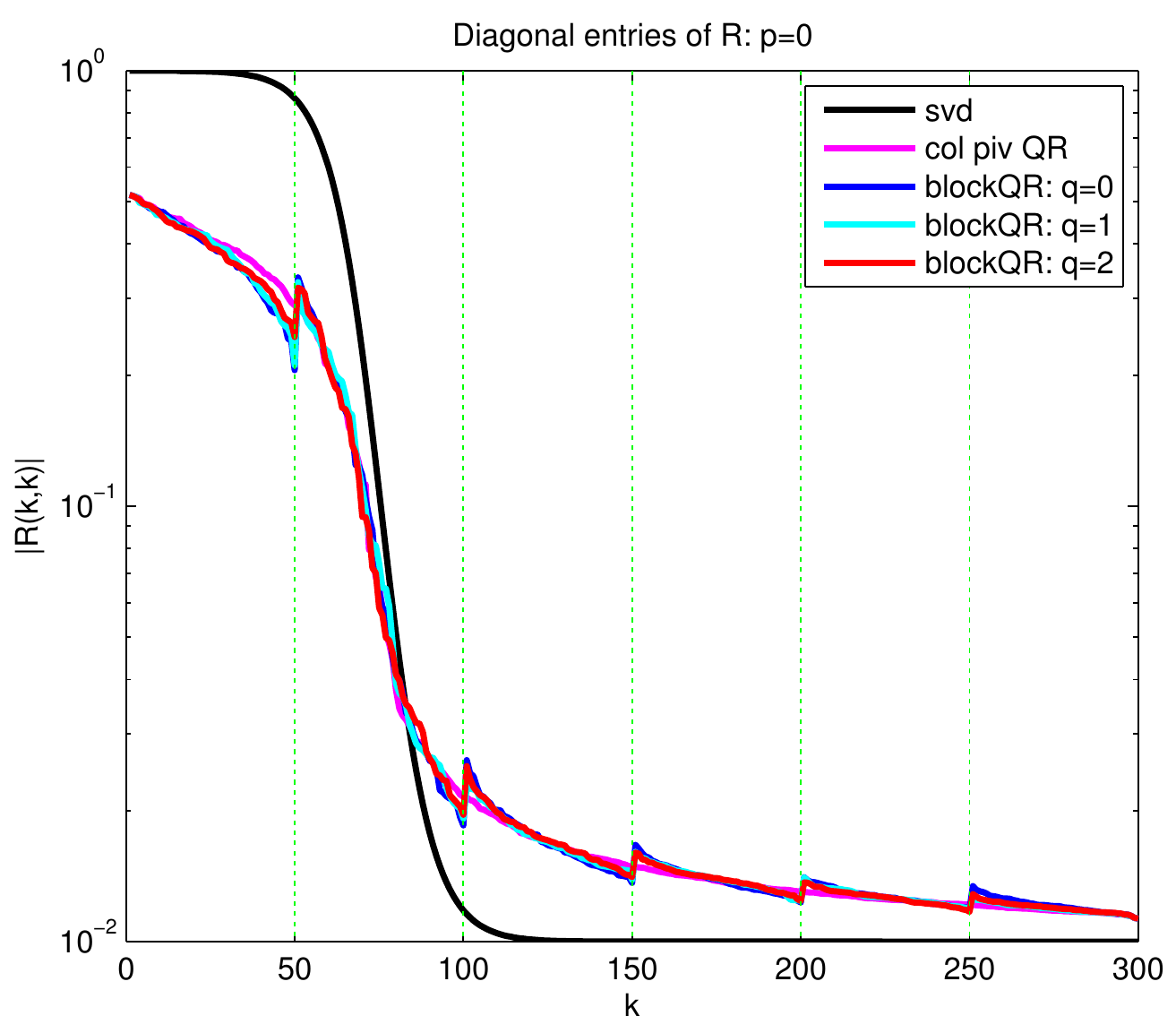}
\caption{Comparison of $|\mtx{R}(k,k)|$ for the matrix $\mtx{R}$ resulting from Method 1,
and the singular values of $\mtx{A}^{\rm S}$. We would ideally like these values to be close.}
\label{fig:diag_method1_matrix3}
\end{figure}

\begin{figure}
\includegraphics[width=0.95\textwidth]{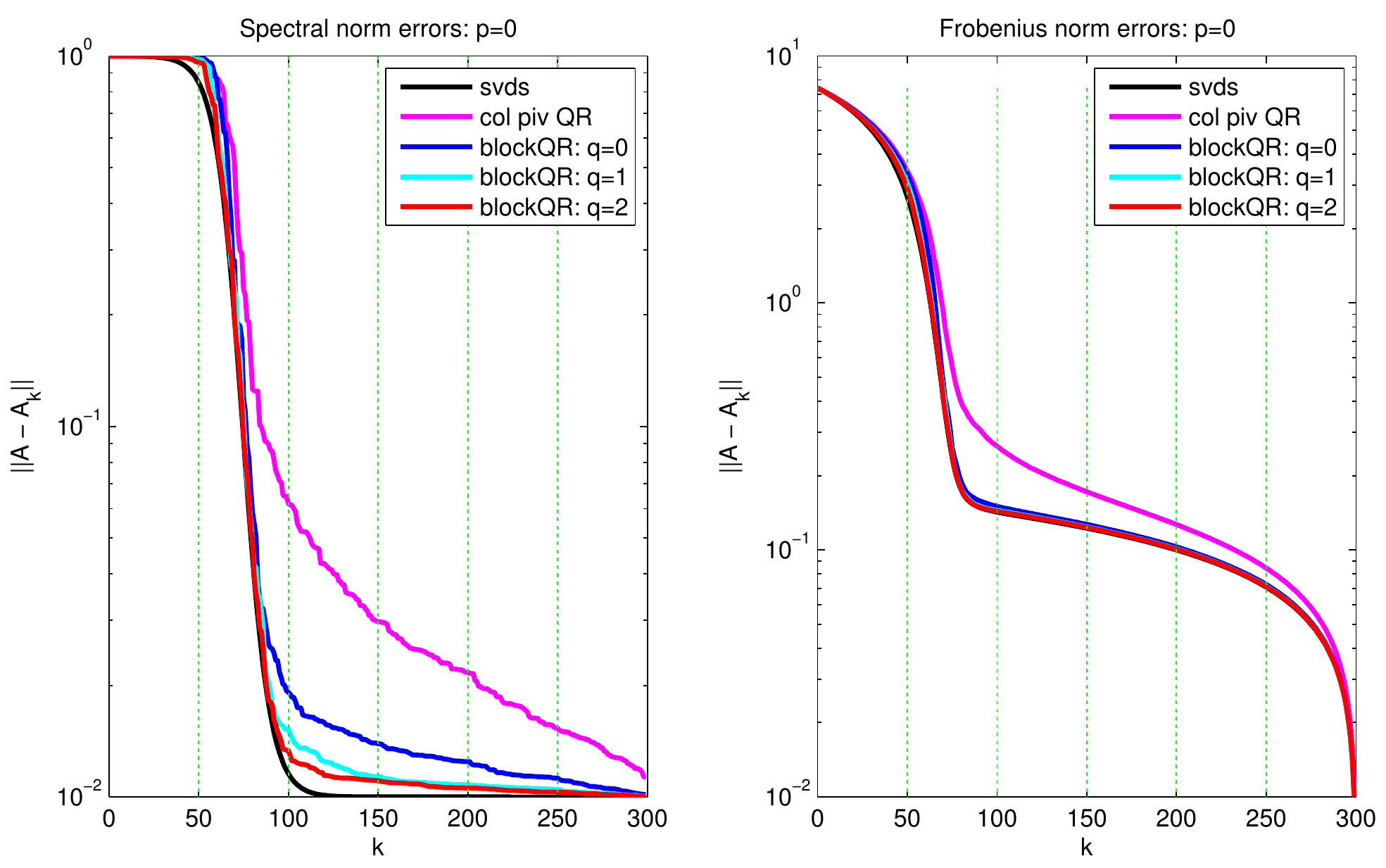}
\caption{Approximation errors from Method 2 (\texttt{blockQR} with a Householder pivoting matrix) for
$\mtx{A} = \mtx{A}^{\rm S}$, cf.~Section \ref{sec:numS}.
Also included are the errors resulting from a truncated SVD (which are theoretically optimal), and
a column pivoted QR factorization.}
\label{fig:err_method2_matrix3}
\end{figure}

\begin{figure}
\includegraphics[width=110mm]{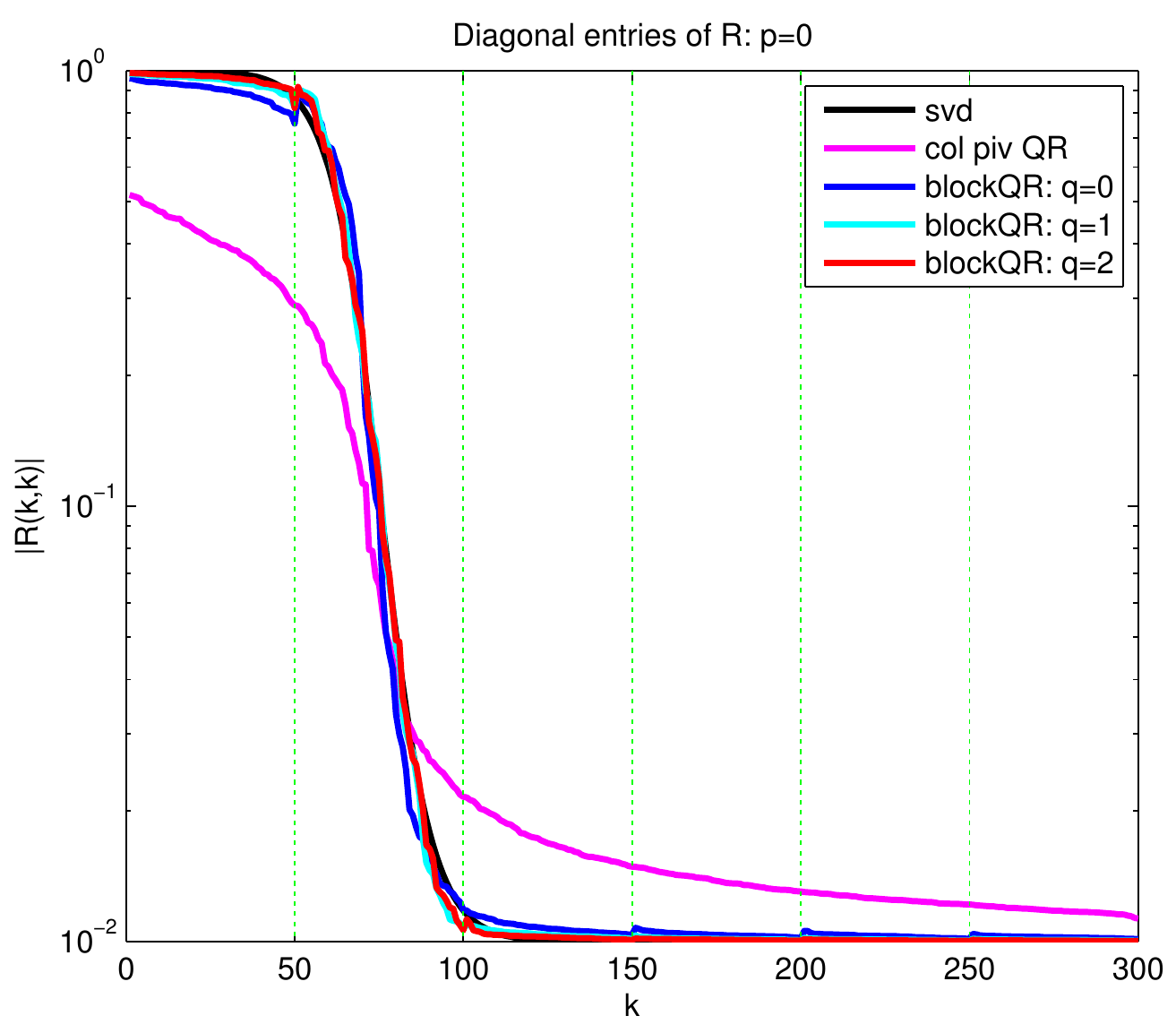}
\caption{Comparison of $|\mtx{R}(k,k)|$ for the matrix $\mtx{R}$ resulting from Method 2,
and the singular values of $\mtx{A}^{\rm S}$.}
\label{fig:diag_method2_matrix3}
\end{figure}

\begin{figure}
\includegraphics[width=0.95\textwidth]{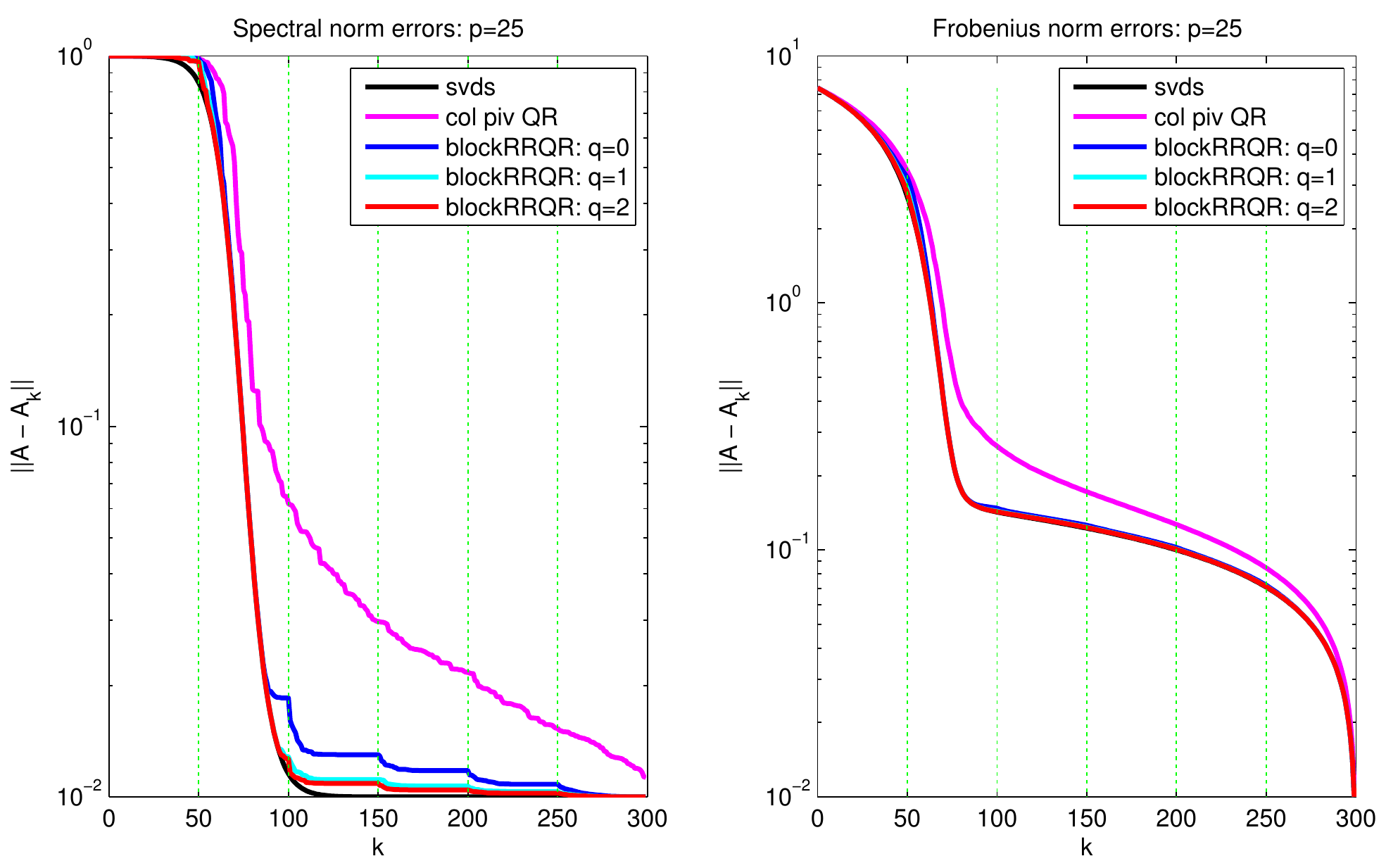}
\caption{Approximation errors from Method 3 (\texttt{blockRRQR} with $p=25$) for
$\mtx{A} = \mtx{A}^{\rm S}$, cf.~Section \ref{sec:numS}.
Also included are the errors resulting from a truncated SVD (which are theoretically optimal), and
a column pivoted QR factorization.}
\label{fig:err_method3_matrix3}
\end{figure}

\begin{figure}
\includegraphics[width=110mm]{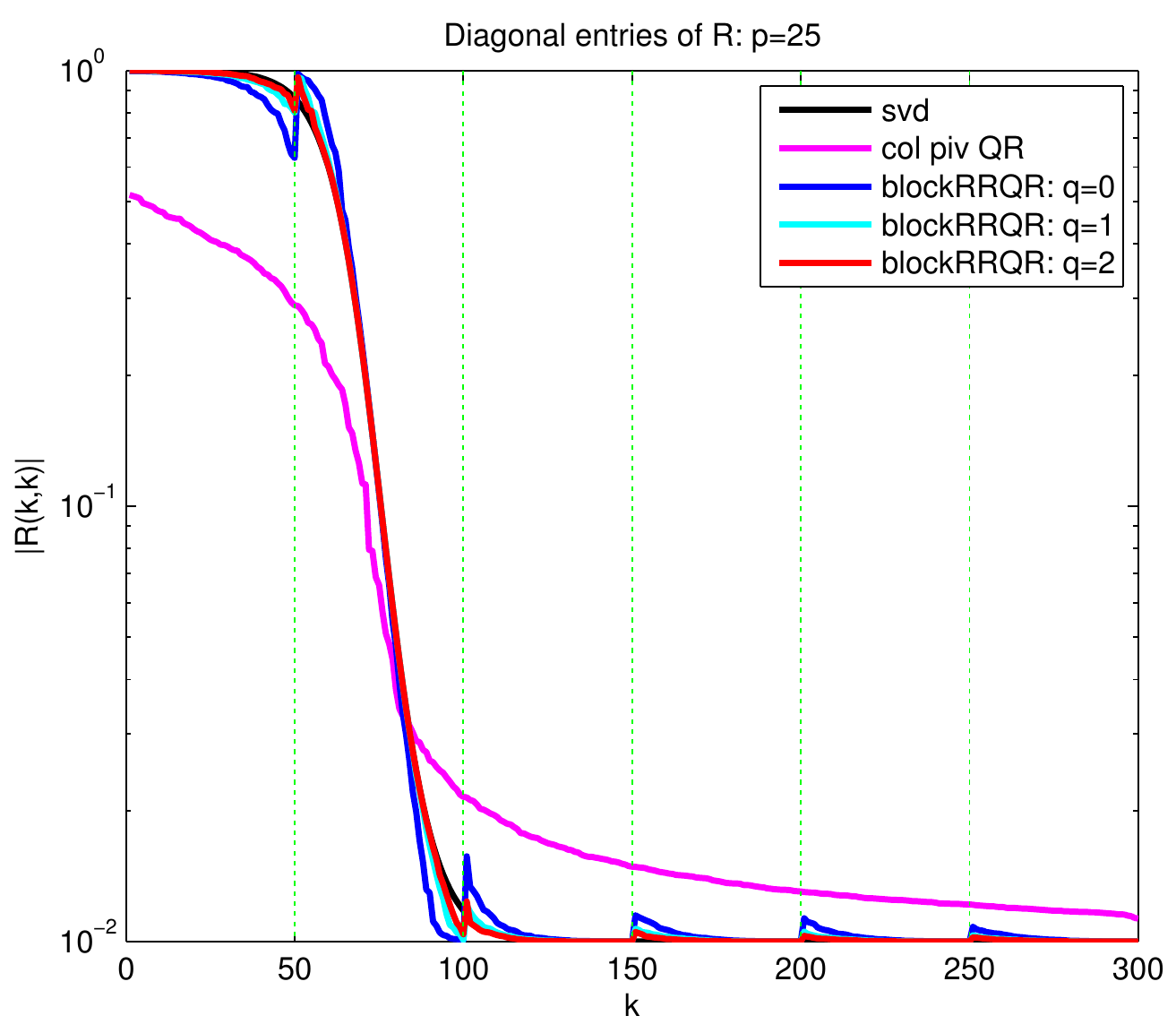}
\caption{Comparison of $|\mtx{R}(k,k)|$ for the matrix $\mtx{R}$ resulting from Method 3,
and the singular values of $\mtx{A}^{\rm S}$.}
\label{fig:diag_method3_matrix3}
\end{figure}

We next repeat all experiments conducted in Section \ref{sec:numfast}, but now for a
square matrix $\mtx{A}^{\rm S}$ given by
$$
\mtx{A}^{\rm S} = \mtx{U}\mtx{D}\mtx{V}^{*},
$$
where $\mtx{U}$ and $\mtx{V}$ are orthonormal (drawn at random from a uniform distribution),
and $\mtx{D}$ is a diagonal matrix whose entries are the singular values of $\mtx{A}^{\rm S}$,
as shown in Figure \ref{fig:err_method1_matrix3}. The singular values of this matrix are chosen
to be flat for a while, then decrease rapidly, and then level out again. In other words, the
``tail'' of the singular values is now heavy and exhibits no decay, which is known to be a
particularly challenging environment for the randomized sampling scheme.
The matrix is again of size $300\times 300$,
and we used a block size of $b=50$. The results are shown in Figures \ref{fig:err_method1_matrix3} -- \ref{fig:diag_method3_matrix3}.

The results for this example are qualitatively very similar as what we saw in Sections \ref{sec:numfast}
and \ref{sec:numgauss}. This examples illustrates particularly well that the randomized schemes
are much better at approximating matrices in the Frobenius norm than in the spectral norm. Moreover,
in this example, we see a strong improvement in going from a permutation matrix as the pivoting matrix
to Householder pivoting matrices.

\clearpage

\clearpage

\section{Conclusions}

We have described techniques for efficiently computing a QR factorization
$\mtx{A}\mtx{P} = \mtx{Q}\mtx{R}$ of a given matrix $\mtx{A}$. The main innovation
is the use of randomized sampling to determine the ``pivot matrix'' $\mtx{P}$.
The randomization allows us to \textit{block} the factorization which we expect
will substantially accelerate execution speed, in particular in communication
constrained environments such as a matrix processed on a GPU or a distributed
memory parallel machine, or stored out-of-core.

We also discussed a variation of the QR factorization where we allow the matrix
$\mtx{P}$ to be a product of Householder reflectors (as opposed to a permutation matrix
in the classical setting). We demonstrated through numerical experiments that this
generalization leads to dramatic improvements in the approximation error obtained by
truncated factorizations.

\vspace{5mm}

\textit{\textbf{Acknowledgements:}} The research reported was supported by DARPA,
under the contract N66001-13-1-4050, and by the NSF, under the contract DMS-1407340.

\bibliography{main_bib}

\providecommand{\bysame}{\leavevmode\hbox to3em{\hrulefill}\thinspace}
\providecommand{\MR}{\relax\ifhmode\unskip\space\fi MR }
% \MRhref is called by the amsart/book/proc definition of \MR.
\providecommand{\MRhref}[2]{%
  \href{http://www.ams.org/mathscinet-getitem?mr=#1}{#2}
}
\providecommand{\href}[2]{#2}
\begin{thebibliography}{1}

\bibitem{1987_bischof_WY}
Christian Bischof and Charles Van~Loan, \emph{The wy representation for
  products of householder matrices}, SIAM Journal on Scientific and Statistical
  Computing \textbf{8} (1987), no.~1, s2--s13.

\bibitem{1987_chan_RRQR}
Tony~F Chan, \emph{Rank revealing qr factorizations}, Linear Algebra and Its
  Applications \textbf{88} (1987), 67--82.

\bibitem{1936_eckart_young}
Carl Eckart and Gale Young, \emph{The approximation of one matrix by another of
  lower rank}, Psychometrika \textbf{1} (1936), no.~3, 211--218.

\bibitem{golub}
Gene~H. Golub and Charles~F. Van~Loan, \emph{Matrix computations}, third ed.,
  Johns Hopkins Studies in the Mathematical Sciences, Johns Hopkins University
  Press, Baltimore, MD, 1996.

\bibitem{gu1996}
Ming Gu and Stanley~C. Eisenstat, \emph{Efficient algorithms for computing a
  strong rank-revealing {QR} factorization}, SIAM J. Sci. Comput. \textbf{17}
  (1996), no.~4, 848--869. \MR{97h:65053}

\end{thebibliography}
\bibliographystyle{amsplain}

\end{document}